\documentclass[11pt,a4paper,twoside]{article}
\usepackage{a4}
\usepackage{bm, amsmath, amssymb, amsthm} 
\usepackage{graphicx}        

\def\bu{u}

\def\dx{\partial_x\/}
\def\ddx{\partial^2_{xx}\/}

\def\<{\langle}
\def\>{\rangle}
\def\eps{\varepsilon}
\def\NN{\mathbb{N}}
\def\ZZ{\mathbb{Z}}
\def\RR{\mathbb{R}}

\def\calf{\mathcal{F}}

\def\tr{\operatorname{tr}\,}
\def\id{\operatorname{id}}
\def\Div{\operatorname{div}}

\newcommand\Ric{\operatorname{Ric}}
\newcommand\vol{\operatorname{vol}}

\def\eq{\hspace*{-.1mm}&=&\hspace*{-.1mm}}
\def\plus{\hspace*{-.1mm}&+&\hspace*{-.1mm}}

\def\dt{\partial_t}

\newtheorem{definition}{Definition}
\newtheorem{example}{Example}
\newtheorem{remark}{Remark}
\newtheorem{lemma}{Lemma}
\newtheorem{proposition}{Proposition}
\newtheorem{theorem}{Theorem}

\author{Vladimir Rovenski
\thanks{
        Mathematical Deptartment, University of Haifa,
        E-mail: rovenski@math.haifa.ac.il}
\footnote{
The author would like to thank Pawe\l~Walczak (University of \L\'od\'z)
for helpful corrections concerning the manuscript,
and Igor Gaissinski (Technion, Haifa) for discussion of Section~\ref{subsec:parabPDE}.
The work was supported by the Marie-Curie actions grant, EU-FP7-P-2010-RG, No.~276919.}
}

\title{Extrinsic geometric flows on foliated manifolds,~III}

\begin{document}

\date{}

\maketitle


\begin{abstract}
We study the geometry of a codimension-one foliation with a time-dependent Riemannian metric.
The~work begins with formulae concerning deformations of geometric quantities as the Riemannian metric varies
along the leaves of the foliation.
Then the Extrinsic Geometric Flow depending on the second fundamental form of the foliation is introduced.
Under suitable assumptions, this evolution yields the second order parabolic PDEs, for which the
existence/uniquenes and in some cases convergence of a solution are shown.
Applications to the problem of prescribing mean curvature function of a codimension-one foliation,
and examples with harmonic and umbilical foliations (e.g., foliated surfaces) and with twisted product metrics are given.

\vskip1.5mm\noindent
\textbf{Keywords}: foliation; Riemannian metric; second fundamental form; extrinsic geometric flow; mean curvature; umbilical; harmonic;
heat equation; double-twisted product

\vskip1mm\noindent
\textbf{Mathematics Subject Classifications (2010)} Primary 53C12; Secondary 53C44
\end{abstract}

\section{Introduction}
\label{sec:1-0}

A~\textit{Geometric Flow}  (GF) is an evolution of a given geometric structure under a differential equation related to a functional on a manifold, usually associated with some curvature.
The popular GFs in mathematics are the \textit{heat flow}, the \textit{Ricci flow} and the \textit{Mean Curvature flow}. GFs play an essential role in many fields of mathematics and  physics (see \cite{bre}, \cite{ta}, \cite{to} for the heat flow and Ricci flow; and
\cite{LMR}, \cite{rw-m}\,--\,\cite{rw4} for foliated Riemannian manifolds).
They all correspond to dynamical systems in the infinite-dimensional space of all appropriate geometric structures
(e.g., metrics) on a given manifold.
GF equations are quite difficult to solve in all generality, because of their nonlinearity.
Although the short time existence of solutions is guaranteed by the parabolic or hyperbolic nature of the equations, their (long time) convergence to canonical geometric structures is analyzed under various conditions.
A GF on a foliated manifold is called \textit{extrinsic}, if the evolution depends on the second fundamental tensor of the foliation.

 One of the principal problems of extrinsic geometry of foliations reads as follows:
 \textit{Given a foliation $\calf$ on a manifold $M$ and an extrinsic geometric property (P), does there exist a Riemannian metric
$g$ on $M$ such that $\calf$ enjoys (P) w.\,r.\,t.~$g$?}

Such problems (first posed by H.\,Gluck for geodesic foliations) were studied already in the 1970's when D.\,Sullivan
provided a topological condition (called {\it topological tautness}) for a foliation, equivalent to
the existence of a Riemannian metric making all the leaves minimal.

Let $M^{n+1}$ be a connected smooth manifold with a codimension-one foliation $\calf$, and a vector field $N$ transversal to $\calf$.
A Riemannian metric  $g$ on $M$ is \textit{adapted} if $g(N, T\calf)=0$ and $g(N,N)=1$.
Denote by $\mathcal{M}$ the space of adapted Riemannian metrics of finite volume for $(M,\calf,N)$.

 The definition of the $\calf$-\emph{truncated} $(r,k)$-\emph{tensor field} $\hat S$
(where $r=0,1$, and $\,\widehat{}\,$ is the $T\calf$-component) will be helpful,
 $\hat S(X_1,\dots,X_k) = S(\hat X_1,\dots,\hat X_k)\ (X_i\in TM)$.
For example, $\hat g$ is the $\calf$-truncated metric tensor $g$ on $M$,
i.e., $\hat g(T\calf, T\calf)=g(T\calf, T\calf)$, and $\hat g(N,\cdot)=0$.
Another example is the {second fundamental form} $b(X,Y)=g(\nabla_X Y, N)$, where $X,Y\in T\calf$,
and the {Weingarten operator} $A(X) = -(\nabla_X N)^\top$ of $\calf$ w.\,r.\,t. $g$,
extended on $TM$ by $b(N,\cdot)=0$ and $A(N)=0$.
Denote $b_m$ the symmetric $(0,2)$-tensor dual to
$A^m$, hence, $b_m(X,Y)= g(A^m(X),Y)$.

We call the \textit{Extrinsic Geometric Flow} (EGF)  on $(M,\calf,N)$ (see Section~\ref{sec:1-def})
a family  $g_t$ of adapted Riemannian metrics
(i.e., $g_t(N, T\calf)=0$, $g_t(N,N)=1$) satisfying along $\calf$
\begin{equation}\label{eq1-general}
 \dt g_t = S(b_t).
\end{equation}
Here $S(b)$ is a symmetric $\calf$-truncated $(0,2)$-tensor expressed in terms of $b$.
 Recently, the author and P.\,Walczak
studied EGFs on codimension-one foliations, see  \cite{rw3}\,--\,\cite{rw4} and \cite[Parts 2 and 3]{rw-m},
\begin{equation}\label{eq1-a}
 \dt g_t= h(b_t),
\end{equation}
where the symmetric $\calf$-truncated $(0,2)$-tensor $h(b)$ is given by
\[
  h(b)  =\sum\nolimits_{m=0}^{n-1} f_m\,{b}_m,
\]
and $f_m$ are either functions on $M\times\RR$ or symmetric functions of the principal curvatures of $\calf$.
 The~EGFs (\ref{eq1-a}) reduce to the first order PDEs and yield quasi-linear hyperbolic PDEs.
 In the paper we introduce and study~EGFs (\ref{eq1-general}) of the form
\begin{equation}\label{eq1}
 \dt g_t=\nabla^t_N h(b_t).
\end{equation}
The $\calf$-truncated symmetric tensor $h(b)$ can be expressed in terms of the first partial derivatives of the metric~$g$.
Therefore $g\mapsto\nabla_N h(b)$ is a second-order differential operator.
Hence, (\ref{eq1}) reduce to the second order PDEs and yield quasi-linear parabolic PDEs  (as will be shown in Sections~\ref{subsec:2-3-3}--\ref{sec:proofs}).

EGFs preserve the properties of codimension-one foliations:
\textit{umbilical} ($A=\lambda\,\widehat\id$, see Proposition~\ref{P-08b}),
\textit{totally geodesic} ($A=0$), and \textit{Riemannian} ($\nabla_N N=0$).

We show that EGFs serve as a tool for studying the following \textbf{question} on geometry of foliations:

 \emph{Under what conditions on $(f_m)$ and $(M, \calf, g_0)$ EGF metrics $g_t$
 converge to one for which $\calf$ enjoys a given extrinsic geometric property, e.g., is umbilical, totally geodesic, harmonic,~etc?}

\vskip1mm
The structure of the paper is as follows.
Section~\ref{sec:main-res} collects main results of the paper:
in Theorems~\ref{T-01}--\,\ref{T-02} we study local existence and uniqueness of EGFs
(introduced in Section~\ref{sec:1-def});
in Theorems~\ref{C-nablab}--\,\ref{C-tau-sigma-k} and Propositions~\ref{C-twisted}--\ref{C-tauf} we apply the method to the problem of prescribing mean curvature of a foliation.
Section~\ref{subsec:tvar} contains formulae for the deformation of
extrinsic geometric quantities of a foliation as the Riemannian metric varies along $\calf$.
Section~\ref{subsec:2-8-3} contains examples for EGFs on foliated surfaces.
Section~\ref{subsec:2-3-2} (Appendix) collects the necessary facts about parabolic PDEs
and extend results of \cite{rw-m} about the generalized companion matrix.

\section{Main Results}
\label{sec:main-res}

The \textit{power sums} of the principal curvatures $k_1, \ldots, k_n$ of $\calf$
(the eigenvalues of $A$) are given~by
\begin{equation*}
 \tau_j = k_1^j + \ldots + k_n^j = \tr(A^j)\qquad (j\ge0).
\end{equation*}
The $\tau$'s can be expressed using the \textit{elementary symmetric functions}
  $\sigma_j = \sum\nolimits_{i_1 <\ldots <i_j} k_{i_1}\cdot\ldots\cdot k_{i_j}$ where $0\le j\le n$.
Notice that $\sigma_0=1$ and $\tau_0=n$.
Set $\overrightarrow{\!\tau}=(\tau_1,\ldots, \tau_n)$
and $\overrightarrow{\!\sigma}=(\sigma_1,\ldots, \sigma_n)$.

The choice of the right-hand side (in EGF equations) for $h(b)$ seems to be natural; the powers ${b}_m$ are the only $(0,2)$-tensors which can be obtained algebraically from the  second fundamental form $b$, while $\tau_1, \ldots, \tau_n$ (or, equivalently, $\sigma_1, \ldots, \sigma_n$) generate all the scalar invariants of extrinsic geometry.
Powers ${b}_m$ with $m > 1$ in (\ref{eq1}) or (\ref{eq1-Lie}) are meaningful; for example, EGFs
\begin{equation}\label{eq1-NT1}
\dt g_t =\nabla^t_N\,T_r(b_t),
\end{equation}
produced by the Newton transformation $T_r$, see Example~\ref{Ex-NT}, depend on all ${b}_m\ (m\le r)$.

\begin{remark}\rm
Let $\mathcal{L}_N$ denote the Lie derivative w.\,r.\,t. $N$.
Recall that $\mathcal{L}_N f=N(f)$ for a function $f$ on $M$, and $\mathcal{L}_N X=[N,X]$ for a vector field.
If $\alpha$ is any $(0,j)$-tensor, then
\[
 \mathcal{L}_N(X_1,\ldots,X_j)= N(\alpha(X_1,\ldots,X_j)) -\sum\nolimits_j \alpha(X_1,\ldots,X_{i-1}, [N,X_i], X_{i+1},\ldots,X_j)
\]
for any vector fields $X_1,\ldots,X_j$.
The following EGF is similar to (\ref{eq1}):
\begin{equation}\label{eq1-Lie}
 \dt g_t=\mathcal{L}_N h(b_t).
\end{equation}
Since $h(b)$ is a symmetric $\calf$-truncated $(0,2)$-tensor, $\mathcal{L}_N h(b)$ preserves this property:
\[
 \mathcal{L}_N h(b)(N, X)= N(h(b)(N,X)) - h(b)([N,N],X) - h(b)(N,[N,Y]) = 0\quad (X\in TM).
\]
Theorem~\ref{T-01} (and Theorem~\ref{T-02} in what follows) about EGF (\ref{eq1}), are also valid for EGF (\ref{eq1-Lie}).
\end{remark}

 Define the quantities $\mu_{mij}$, depending on the principal curvatures $(k_i)$ of $A$,~by
\[
 \mu_{m;ij}=\sum\nolimits_{\alpha=0}^{m-1} k_i^\alpha k_j^{m-1-\alpha}\qquad (1\le i\le j\le n,\quad 1\le m<n).
\]
Remark that $\mu_{1;ij}=1$ for all $i,j$. The \textit{ellipticity condition} for EGFs takes the form
\begin{equation}\label{E-cond-a}
 \sum\nolimits_{m=1}^{n-1} f_m\,\mu_{m;ij}<0,\qquad (1\le i\le j\le n).
\end{equation}
In case of umbilical foliations (i.e., $A=\lambda\,\hat\id$), the condition (\ref{E-cond-a}) reads as
 $\sum\limits_{m=1}^{n-1} m\,f_m\,\lambda^{m-1}<0$.

\vskip1.5mm
Theorem~\ref{T-01} concerns the local existence and uniqueness of EGFs with $f_m$ of type (a).

\begin{theorem}\label{T-01}
Let $(M, g_0)$ be a closed Riemannian manifold with a codimen\-sion-one foliation $\calf$, and let $f_m\in C^{2}(M\times\RR)$. If (\ref{E-cond-a}) is satisfied on $M$ for $t=0$, then there is a unique smooth solution $g_t$ of (\ref{eq1}) defined on a positive time interval $[0,\eps)$.
\end{theorem}

So, given $f_m\in C^{2}(\RR)$, if (\ref{E-cond-a}) holds for $t=0$ then there is a unique smooth local solution
to
\[
 \dt g_t=\sum\nolimits_{m=0}^{n-1} f_m(t)\nabla^t_N\,b_m.
\]

The \textit{generalized companion matrix} $B_{n,1}$ is defined in \cite{rw-m}, see Section~\ref{subsec:GCmatrix}.
Consider the matrices
\begin{equation}\label{E-BBB-main}
 \widetilde{B}_n=\sum\nolimits_{m=1}^{n-1} \frac m2\,f_m(\overrightarrow{\!\tau})\,(B_{n,1})^{m-1},\qquad
\tilde A = (\tilde{A}_{ij}),\quad\tilde A_{ij}=\frac i2\sum\nolimits_{m=0}^{n-1} \tau_{i+m-1}\,f_{m,\tau_j}(\overrightarrow{\!\tau}).
\end{equation}

The next Theorem~\ref{T-02} about EGFs with $f_m$ of type~(b) is based on Theorem~\ref{T-01}.

\begin{theorem}[\textbf{Short time existence}]\label{T-02}
Let $(M, g_0)$ be a closed Riemannian manifold with a codimensi\-on-one foliation $\calf$
and a unit normal $N$. Given  $f_m\in C^{2}(\RR^{n})$, suppose that  at $t=0$ the condition (\ref{E-cond-a}) is satisfied
and the matrix $\widetilde{B}_n+\widetilde{A}$, see (\ref{E-BBB-main}), is negative definite on $M$.
Then there is a unique smooth solution $g_t$ of (\ref{eq1}) defined on a positive $t$-interval $[0,\eps)$.
\end{theorem}

\begin{example}\label{ex-1-Th2}\rm

 (i) Let $f(\overrightarrow{\!\tau})\ge c>0$. Consider EGF
\begin{equation*}
 \dt g_t = -f(\overrightarrow{\!\tau})\,\nabla^t_N\,b_t.
\end{equation*}
The corresponding system for $\tau$'s (see Lemma~\ref{L-btAt2} in Section~\ref{subsec:tvar}
with $S=f(\overrightarrow{\!\tau})\,\nabla_N\,b$) is
\[
 \dt\tau_i =\frac{1}{2}\,f(\overrightarrow{\!\tau})\,N(N(\tau_{i}))+a_i,\qquad i>0,
\]
where
$a_i=\frac12\,N(f(\overrightarrow{\!\tau}))N(\tau_i)
-\frac{i}2 f(\overrightarrow{\!\tau})\tr(\nabla_N(A^{i-1})\nabla_N(A))$.
The above in-homogeneous parabo\-lic system for $\tau$'s along $N$-curves (see Appendix),
has a unique solution for $0\le t\le\eps$.

\vskip1.5mm
 (ii) Consider EGF with $f_m$ of type (b), given along $\calf$ by
\begin{equation*}
 \dt g_t = \nabla^t_N\big(f(\overrightarrow{\!\tau})\,b_t\big),
\end{equation*}
where
$\frac{1}{2}\,f(\overrightarrow{\!\tau})\,\hat\id+\widetilde{A}<0$
($\widetilde{A}_{ij}=\frac i2\,\tau_i\,f_{, j}$)
for all $\overrightarrow{\!\tau}(q)$ and $a\in M$.
The system
for $\tau$'s (see Proposition~\ref{P-tauAk-a} in Section~\ref{P-tauAk-a}) is parabolic quasi-linear
\[
 \dt\tau_i=-\big(\frac{1}{2}\,f(\overrightarrow{\!\tau})\,\hat\id+\widetilde{A}(\overrightarrow{\!\tau})\big)N(N(\tau_{i}))
 +a_i,\qquad i>0,
\]
where
 $a_i=\frac{i}2\,f(\overrightarrow{\!\tau})\tr(\nabla^t_N(A^{i-1})\nabla^t_N(A))
 -N(f(\overrightarrow{\!\tau}))\,N(\tau_{i})
 -\frac{i}2\,\tau_{i}\!\sum\nolimits_{a,b=1}^{n}\!f_{,ab}(\overrightarrow{\!\tau}) N(\tau_a)N(\tau_b)$.
The solution is unique and exists as long as $\|\overrightarrow{\!\tau}\|$ is bounded.
\end{example}

A foliation with vanishing $\tau_1$ is called \textit{harmonic}; its leaves are minimal submanifolds.
A foliation $\calf$ is \textit{taut} if there is at least one metric on $M$ for which $\calf$ is harmonic.
The known proofs of the existence of such metrics use the Hahn--Banach Theorem and are not constructive.
 In particular, \textit{$\calf$ is taut if and only if there is a closed transversal curve that intersects each leaf}, see~\cite{cc1}, \cite{rsw}.
For example, a Reeb foliation on $S^3$ is not taut.

Theorem~\ref{C-nablab} (based on Theorem~\ref{T-01}) helps us to find in some cases
a one-parameter family of metrics converging to the metric for which $\calf$ is harmonic.

\begin{theorem}\label{C-nablab}
The following EGF, see (\ref{eq1-Lie}), on a closed Riemannian manifold $(M,g_0)$:
\begin{equation}\label{eq1-1}
 \dt g_t = -2\,\mathcal{L}_N\,b_t
\end{equation}
admits a unique smooth global solution $g_t\ (t\ge0)$. Moreover, if
$N$-curves compose a fibration $S^1\overset{i}\hookrightarrow M\overset{\pi}\to B$,
then the metrics $g_t$ approach to a smooth metric $\bar g$ as $t\to\infty$, and $\calf$ is $\bar g$-harmonic.
\end{theorem}

\begin{remark}\label{R-replace}\rm
 One may replace (\ref{eq1-1}) in Theorem~\ref{C-nablab} by the PDE (\ref{eq1-1a}) in what follows.
 By Example~\ref{ex-Reeb1} (in Section~\ref{subsec:2-8-3}), the claim is false without condition
``$N$-curves compose a fibration".
\end{remark}

Let $M=M_1\times S^1$ be the product of a Riemannian manifold $(M_1,g_1)$ and a circle $(S^1,g_2)$.
For a positive function $f\in C^1(M)$, the \textit{metric of a twisted product} $M_1\times_{f}\!S^1$ is given by
$g=(f^2 g_1)\oplus g_2$. The foliation $\calf=M_1\times\{y\}$ is umbilical with the
\textit{mean curvature vector} $H=-\nabla(\log f)^\perp$.

\begin{proposition}\label{C-twisted}
Let $M_1\times_{f}\!S^1$ be a twisted product of a closed Riemannian manifold $(M_1, g_1)$ and a circle $(S^1,g_2)$.
Then EGF (\ref{eq1-1}) admits a unique smooth solution $g_t\in\mathcal{M}$ for all $t\ge0$,
consisting of twisted product metrics on  $M_1\times_{f_t}\!S^1$.
As $t\to\infty$, the metric $g_t$ converges to the metric $\bar g$ of the product
$(M_1,\bar f_1\cdot g_1)\times(S^1,g_2)$, where $\bar f(x)=\int_{S^1} f(0, x, y)\,d y_g $.
\end{proposition}

\begin{example}\rm
Given $f:T^2\to\RR_+$, let $S^1\times_{f}\!S^1$ be a twisted product of two circles $(S^1, g_i)$.
Then EGF (\ref{eq1-1}) has a unique smooth solution $g_t\in\mathcal{M}$ for all $t\ge0$,
consisting of twisted product metrics on the torus~$T^2$.
As $t\to\infty$, the metric $g_t$ converges in $C^0$
to the flat metric $\bar g$ on $T^2$.
\end{example}

\begin{proposition}\label{C-ftau}
Given $f\in C^\infty(\RR^n)$, suppose that EGF on $M$, either (i) or (ii),
\begin{equation}\label{E-FPsi}
 (i)~\dt g_t=-N(f(\overrightarrow{\!\tau}))\,\hat g_t,\qquad
 (ii)~\dt g_t=-N(f(\overrightarrow{\!\sigma}))\,\hat g_t
\end{equation}
satisfies the inequality for $g_0$ (where functions $F_{k}, \Psi_{k}$ are defined in Proposition~\ref{P-conf}):
\[
 (i)~\sum\nolimits_{k=1}^n k\,F_{k-1}(\tau_1) f_{,k}(\overrightarrow{\!\tau})>0, \quad
 (ii)~\sum\nolimits_{k=1}^n (n-k+1)\,\Psi_{k-1}(\sigma_1) f_{,k}(\overrightarrow{\!\sigma})>0.
\]
Then (\ref{E-FPsi}) admits a unique smooth solution $g_t$ defined on a positive $t$-interval $[0,\eps)$.
\end{proposition}

\begin{example}\label{ex-1-4a}\rm
Notice that for $f=\frac2n\,\tau_1$, EGF (\ref{E-FPsi}) coincides with
\begin{equation}\label{eq1-1a}
 \dt g_t = -\frac2n\,N(\tau_1)\,\hat g_t.
\end{equation}
 We~ shall illustrate Proposition~\ref{C-ftau} for $f=\frac2n\,\tau_k\ (k>1)$ and $f=\frac2n\,\sigma_k\ (k>1)$.
 Consider EGFs
\begin{equation}\label{eq1-1ak}
 {\rm  (i)}~\dt g = -\frac2n\,N(\tau_k)\,\hat g,\qquad
 {\rm (ii)}~\dt g = -\frac2n\,N(\sigma_k)\,\hat g.
\end{equation}
If
(i) $F_{k-1}(\tau_1)>0$, (ii) $\Psi_{k-1}(\sigma_1)>0$,
then the corresponding PDE is parabolic:
\begin{equation*}
 {\rm (i)}~\dt\tau_1 = N\big(F_{k-1}(\tau_1)\,N(\tau_1)\big),\qquad
 {\rm (ii)}~\dt\sigma_1 = \frac{n-k+1}{n}\,N\big(\Psi_{k-1}(\sigma_1)\,N(\sigma_1)\big).
\end{equation*}
Solving these systems and using relations $\tau_k=F_{k}(\tau_1)$ and $\sigma_k=F_{k}(\sigma_1)$, see Proposition~\ref{P-conf}
in Section~\ref{subsec:tvar}, one may find a local solution of EGFs (\ref{eq1-1ak}).
\end{example}

Using Propositions~\ref{C-ftau} and \ref{P-heat-eq} (in Section~\ref{subsec:parabPDE}), one may show the following.

\begin{theorem}\label{C-tau-sigma-k}
If for some real $c>0$ the following ellipticity condition holds:
\[
 (i)~c^{-1}\le\sum\nolimits_{k=1}^n \!k\,F_{k-1}(\tau_1) f_{,k}(\overrightarrow{\!\tau})\le c, \quad
 (ii)~c^{-1}\le\sum\nolimits_{k=1}^n (n-k+1)\,\Psi_{k-1}(\sigma_1) f_{,k}(\overrightarrow{\!\sigma})\le c,
\]
then (\ref{E-FPsi}) admits a unique smooth global solution $g_t\ (t\ge0)$.
Moreover, if $N$-curves compose a fibration $S^1\overset{i}\hookrightarrow M\overset{\pi}\to B$,
then $g_t$ approach to a smooth metric $\bar g$ as $t\to\infty$, and $\calf$ is $\bar g$-harmonic.
\end{theorem}

A compact saturated domain $D\subset M$ is called a $(+)$-\textit{fcd},
if the transverse orientation of $\calf$ is outward everywhere on $\partial D$,
and $(-)$-\textit{fcd} if the transverse orientation of $\calf$ is inward on $\partial D$.
P.~Walczak and G.~Oshikiri studied the question for codimension-one foliations (see \cite{o1997}, \cite{w1984} etc):

\textit{Which function $F$ on $M$ is admissible, i.e., can be the mean curvature of $\calf$ w.\,r.\,t. some Riemannian metric on~$M$?}
It was proven that $F$ is admissible if and only if

\textit{$F(x)>0$ somewhere in any minimal $(+)$-fcd and
$F(y)<0$ somewhere in any minimal $(-)$-fcd}.

\noindent
In view of the integral formula by Reeb, see \cite{rw-m},
\begin{equation}\label{E-IF-Reeb}
 \int_M \tau_1\,d\vol=0,
\end{equation}
the condition $\int_M F\,d\,V=0$ (for some oriented volume form $dV$ on $M$) is necessary.

\vskip1.5mm
In the next proposition we show how to produce  in some cases a one-parameter family of metrics
whose mean curvature (of~$\calf$) converges to a prescribed function.

\begin{proposition}\label{C-tauf}
Given $F\in C^2(M)$ with the property $\int_M F\,d\vol_g=0$ on a closed Riemannian manifold manifold $(M,g)$, the PDE
\begin{equation}\label{eq1-1f}
 \dt g_t = -\frac2n\,N(\tau_1-F)\,\hat g_t
\end{equation}
has a unique smooth
solution $g_t\,(t\ge0)$. Moreover, if
each $N$-curve is dense in $M$,
then  $\lim\limits_{t\to\infty} \tau_1^t=F$.
\end{proposition}

\section{$\calf$-truncated variations of geometric quantities}
\label{subsec:tvar}

Let $\mathcal{M}_1\subset\mathcal{M}$ be the subspace of metrics of unit volume,~and
\[
 \pi:\mathcal{M}\to \mathcal{M}_1,\qquad
 \pi(g)=\bar g=\big(\vol(M,g)^{-2/n}\,\hat g\big)\oplus g^\perp
\]
the $\calf$-conformal projection.
 Let $g_t\in\mathcal{M}$ (with $0\le t<\eps$) be a family of metrics with $\vol(M,g_0)=1$.
 Set $S_t=\dt g_t$.
Observe that $[X,Y]^\perp$ does not depend on $t$.
The volume form of $g_t$ evolves as, see~\cite{to},
\begin{equation}\label{E-dtdvol}
 \frac{d}{dt}\vol_t=\frac12\,(\tr S^\sharp)\vol_t.
\end{equation}
If $S_t=s_t\,\hat g_t$, where $s_t:M\to\RR$ (i.e., $g_t$ are $\calf$-conformal
and $\hat g_t$ is the $D$-truncated metric $g_t$) then
 $\frac{d}{dt}\vol_t=\frac n2\,s_t\vol_t$.
Hence, the metrics $\tilde g_t=(\phi_t\hat g_t)\oplus g^\perp_t$
with dilating factors $\phi_t=\vol(M,g_t)^{-2/n}$, belong to $\mathcal{M}_1$.

 Recall that the \textit{Levi-Civita connection} $\nabla^t$ of a Riemannian metric $g_t$ on $M$ is given by
\begin{eqnarray}\label{eqlevicivita}
\nonumber
 2\,g_t(\nabla^t_X Y, Z) \eq X(g_t(Y,Z)) + Y(g_t(X,Z)) - Z(g_t(X,Y)) \\
 \plus g_t([X, Y], Z) - g_t([X, Z], Y) - g_t([Y, Z], X)
\end{eqnarray}
for all vector fields $X, Y$ and $Z$ on $M$.
 Since the difference of two connections is always a tensor, $\Pi_t:=\dt\nabla^t$ is a $(1,2)$-tensor field on $(M, g_t)$.
Differentiation (\ref{eqlevicivita}) w.\,r.\,t. $t$ yields the formula, see~\cite{to},
\begin{equation}\label{eq2}
 g_t(\Pi_t(X, Y), Z)=\frac12\big[(\nabla^t_X S_t)(Y,Z)+(\nabla^t_Y S_t)(X,Z)-(\nabla^t_Z S_t)(X,Y)\big]
\end{equation}
for all $X,Y,Z\in\Gamma(TM)$.
If the vector fields $X=X(t),\,Y=Y(t)$ are $t$-dependent, then
\begin{equation}\label{eq2time}
 \dt\nabla^t_X Y = \Pi_t(X, Y) + \nabla^t_X(\dt Y)+ \nabla^t_{\dt X} Y.
\end{equation}
Notice the symmetry $\Pi_t(X, Y)=\Pi_t(Y, X)$ of the tensor $\Pi_t$.

Next, we find the variational formula for $A$, and apply it to
symmetric functions $\tau_j$ and $\sigma_j$ of~$A$.

\begin{lemma}[see \cite{rw-m}]\label{L-btAt2}
Let $g_t\in\mathcal{M}$
and $S=\dt g_{t}$ -- a $\calf$-truncated tensor. Then the Weingarten operator $A$ of $\calf$
(w.\,r.\,t. $N$) and the symmetric functions $\tau_i$ and $\sigma_i$ of $A$ evolve~by
\begin{eqnarray}\label{E-btAt}
 \dt A \eq \frac12\big([A, \, S^\sharp] -\nabla^t_{N} S^\sharp\big),\\
\label{E-At-3}
 \dt\tau_i\eq -\frac i2 \tr(A^{i-1}\nabla^t_{N} S^\sharp),\quad
 \dt\sigma_i= -\frac12\tr(T_{i-1}(A)\nabla^t_{N} S^\sharp),\quad i>0.
\end{eqnarray}
For $S=s\,\hat g$, where $s\in C^1(M)$, we have
\begin{equation*}
 \dt A = -\frac12 N(s)\,\hat{\id} ,\quad
 \dt\tau_i = -\frac i2\,\tau_{i-1} N(s), \quad
 \dt\sigma_i =-\frac12(n-i+1)\,\sigma_{i-1} N(s),\quad i>0.
\end{equation*}
\end{lemma}

\begin{remark}\label{L-btAt}\rm
Let $(M,\,g=\hat g\oplus g^\perp)$ be a Riemannian manifold with a codimen\-sion-one foliation $\calf$ and a unit normal $N$.
Given $\varphi\in C^1(M)$, define a metric $\tilde g = (e^{\,2\/\varphi}\hat g)\oplus g^\perp$.
Then the  second fundamental forms and Weingarten operators of $\calf$ w.\,r.\,t. $\tilde g$ and $g$
are related~by, see \cite{rw-m},
\begin{equation*}
 \tilde b = e^{\,2\/\phi}(b -N(\varphi)\,\hat g),\qquad
 \tilde A = A-N(\varphi)\,\widehat\id.
\end{equation*}
So, if $N(\varphi)=0$ (e.g., if $\varphi$ is constant), then
 (i)~$\tilde b=e^{\,2\/\varphi}\,b$; \
 (ii)~$\tilde A = A$; \
 (iii)~ $\tilde\tau_j=\tau_j$.
 \end{remark}

\begin{proposition}\label{P-conf}
 Let $g_t\in\mathcal{M}$ and $\dt g_{t}=s_t\,\hat g_t$. Then
\begin{eqnarray}\label{E-ex-tau-k}
 &&\hskip-7mm\tau_k^t = F_k(\tau_1^t),
 \quad{\rm where}\quad
 F_k(\tau_1)=\frac1{n^{k-1}}\,\tau_1^k +\sum\nolimits_{i=2}^k
 \frac{k!}{i!(k-i)!}\cdot\frac{\phi_{i}}{n^{k-i}}\,\tau_1^{k-i},\\
 \label{E-ex-sigma-k}
  &&\hskip-8mm\sigma_k^t = \Psi_k(\sigma_1^t),
  \ {\rm where}\
 \Psi_k(\sigma_1)=\frac{(n-1)!}{k!(n-k)!}\cdot\frac1{n^{k-1}}\,\sigma_1^k +\sum\nolimits_{i=2}^k
 \frac{(n-i)!}{(k-i)!(n-k)!}\cdot\frac{\psi_{i}}{n^{k-i}}\,\sigma_1^{k-i}
\end{eqnarray}
 for $k>1$. For small $k$, we have
\begin{eqnarray*}
 F_2(\tau_1)=\frac1n\,\tau_1^2+\phi_{2},\
 F_3(\tau_1)=\frac1{n^2}\,\tau_1^3\!+\frac3n\,\phi_{2}\,\tau_1+\phi_{3},\
 F_4(\tau_1)=\frac1{n^3}\,\tau_1^4\!+\frac6{n^2}\,\phi_{2}\,\tau_1^2\!+\frac4n\,\phi_{3}\,\tau_1+\phi_{4},\ etc.\\
 \Psi_2(\sigma_1)=\frac{n-1}{2\,n}\,\sigma_1^2+\psi_{2},\quad
 \Psi_3(\sigma_1)=\frac{(n-1)(n-2)}{6\,n^2}\,\sigma_1^3\!+\frac{n-2}n\,\psi_{2}\,\sigma_1+\psi_{3},\quad {etc}.
\end{eqnarray*}
The $t$-independent functions $\phi_{k}=F_k(0)$ and $\psi_{k}=\Psi_k(0)$ on $M$ are determined at $t=0$ by recursion.
\end{proposition}

\noindent\textbf{Proof}.
(i) We will prove (\ref{E-ex-tau-k}) by induction.
Notice that the last term (with $i=k$) in (\ref{E-ex-tau-k}) is $\phi_{k}$.
By~Lemma~\ref{L-btAt2}, the system for $\tau$' reads as
\begin{equation}\label{eq1-1tau1-k}
 \dt\tau_1=-\frac n2\,N(s),\quad \dt\tau_2=-\tau_1 N(s),\quad\ldots,\quad
 \dt\tau_k=-\frac k2\,\tau_{k-1}\,N(s),\quad\ldots
\end{equation}
By (\ref{eq1-1tau1-k})$_1$ we have $\dt(\tau_1^2)=-n\,\tau_1\,N(s)$.
Substituting $N(s)$ into (\ref{eq1-1tau1-k})$_2$, we find $\dt(\tau_2-\frac1n\,\tau_1^2)=0$.
Hence, $\tau_2-\frac1n\,\tau_1^2=\phi_{2}$ is $t$-independent, and $F_2(\tau_1)=\frac1n\,\tau_1^2+\phi_{2}$,
see (\ref{E-ex-tau-k}) for $k=2$.

Assume that (\ref{E-ex-tau-k}) holds for some $k\ge2$. Using (\ref{eq1-1tau1-k})$_3$, we find
\begin{equation}\label{E-ex-tau-k1}
 \dt\tau_{k+1}
 =-\frac{k+1}2\,F_k(\tau_1)\,N(s)
 =-\frac{k+1}2\,\Big[\frac1{n^{k-1}}\,\tau_1^k +\sum\nolimits_{i=2}^k \frac{1}{n^{k-i}}\cdot
 \frac{k!}{i!(k-i)!}\,\phi_{i}\,\tau_1^{k-i}\Big]\,N(s).
\end{equation}
By (\ref{eq1-1tau1-k})$_1$ we have $\dt(\tau_1^{j+1})=-\frac{n(j+1)}{2}\,\tau_1^{j}\,N(s)$.
Substituting this into (\ref{E-ex-tau-k1}), we obtain
$\dt\tau_{k+1}=\dt F_{k+1}$, and hence (\ref{E-ex-tau-k}) for $k+1$.

(ii)  We will prove (\ref{E-ex-sigma-k}) by induction.
Notice that the last term (with $i=k$) in (\ref{E-ex-sigma-k}) is $\psi_{k}$.
By~Lemma~\ref{L-btAt2}, the system for $\sigma$' reads as
\begin{equation}\label{eq1-1sigma1-k}
 \dt\sigma_1=-\frac n2\,N(s),\quad
 \dt\sigma_2=-\frac{n-1}{2}\,\sigma_1 N(s),\quad\ldots,\quad
 \dt\sigma_k=-\frac{n-k+1}{2}\,\sigma_{k-1}\,N(s),\quad\ldots
\end{equation}
By (\ref{eq1-1sigma1-k})$_1$ we have $\dt(\sigma_1^2)=-n\,\sigma_1\,N(s)$.
Substituting $N(s)$ into (\ref{eq1-1sigma1-k})$_2$, we find $\dt(\sigma_2-\frac{n-1}{2\,n}\,\sigma_1^2)=0$.
Hence, $\sigma_2-\frac{n-1}{2\,n}\,\sigma_1^2=\psi_{2}$ is $t$-independent, and $\Psi_2(\sigma_1)=\frac{n-1}{2\,n}\,\sigma_1^2+\psi_{2}$, see (\ref{E-ex-sigma-k}) for $k=2$.

Assume that (\ref{E-ex-sigma-k}) holds for some $k\ge2$. Using (\ref{eq1-1sigma1-k})$_3$, we find
\begin{equation}\label{E-ex-sigma-k1}
 \dt\sigma_{k+1}=-\frac{n-k}2\,\Big[\frac{(n-1)!}{k!(n-k)!}\cdot\frac1{n^{k-1}}\,\sigma_1^k
 +\sum\nolimits_{i=1}^k \frac{1}{n^{k-i}}\cdot\frac{(n-i)!}{(k-i)!(n-k)!}\,\psi_{i}\,\sigma_1^{k-i}\Big]\,N(s).
\end{equation}
By (\ref{eq1-1sigma1-k})$_1$ we have $\dt(\sigma_1^{j+1})=-\frac{n(j+1)}{2}\,\sigma_1^{j}\,N(s)$.
Substituting this into (\ref{E-ex-sigma-k1}), we obtain
$\dt\sigma_{k+1}=\dt \Psi_{k+1}$, and hence (\ref{E-ex-sigma-k}) for $k+1$.\qed

\begin{remark}\label{R-FPsi}\rm
We may set $F_0=n,\ F_1(\tau_1)=\tau_1$ and $\phi_1=0$.
Differentiating the relation $\tau_k=F_k(\tau_1)$ with $k\ge1$, see Proposition~\ref{P-conf},
and replacing $\dt\tau_1=-\frac n2\,N(s)$ by Lemma~\ref{L-btAt2}, we find
\begin{equation}\label{eq1-1s-NN}
 \dt\tau_k=F_k'(\tau_1)\,\dt\tau_1=-\frac n2\,F_k'(\tau_1) N(s).
\end{equation}
Comparing (\ref{eq1-1s-NN}) with $\dt\tau_k=-\frac k2\,\tau_{k-1}\,N(s)=-\frac k2\,F_{k-1}(\tau_1)\,N(s)$,
see Lemma~\ref{L-btAt2}, we conclude that
\[
 F_k'(\tau_1)=\frac kn\,F_{k-1}(\tau_1)\quad\Rightarrow\quad
 F_k(\tau_1)=\frac kn\int_0^{\,\tau_1} F_{k-1}(x)\,d\,x +\phi_k.
\]
Similarly, for functions $\Psi_k(\sigma_1)$ in (\ref{E-ex-sigma-k})
we set $\Psi_0=1,\ \Psi_1(\sigma_1)=\sigma_1$  and $\psi_1=0$, and obtain
\[
 \Psi_k'(\sigma_1)=\frac{n-k+1}n\,\Psi_{k-1}(\sigma_1)\quad\Rightarrow\quad
 \Psi_k(\sigma_1)=\frac{n-k+1}n\int_0^{\,\sigma_1} \Psi_{k-1}(x)\,d\,x +\psi_k.
\]
\end{remark}

\begin{lemma}[see \cite{rw-m}]\label{L-nablaNN}
The vector field $Z=\nabla^t_{N} N$ is evolved by $g_t\in\mathcal{M}$ with $S=\dt g_t$ as
\begin{equation}\label{E-nablaNNt2}
  (i)~\dt Z = -S^\sharp(Z),\qquad
 (ii)~\dt Z =-s\,Z\quad {\rm for}\quad S=s\,\hat g.
\end{equation}
In particular, all variations $g_t\in\mathcal{M}$ preserve Riemannian foliations.
\end{lemma}

We will use the following condition for convergence of evolving conformal metrics
(for a general formulation see \cite[Appendix~A]{bre}).

\begin{proposition}\label{P-converge}
 Let $\dt g_t=s_t\,\hat g_t\ (t\ge0)$ be a one-parameter family of $\calf$-conformal Riemannian metrics on a closed manifold $M$ with a codimension-one foliation $\calf$ and a unit normal $N$.
 Define a function $v(t)=\sup_M |s_t|_{g(t)}$ and assume that $\int_0^\infty v(t)\,dt<\infty$.
 Then, as $t\to\infty$, the metrics $g_t$ converge in $C^0$ to a smooth limit Riemannian metric $g_\infty$.
\end{proposition}

\noindent\textbf{Proof}. Our assumptions ensure that the metrics $g_t$ converge in $C^0$ to a symmetric $(0,2)$-tensor $g_\infty$.
Notice that the metrics are uniformly equivalent: $c^{-1}\hat g_0\le \hat g_t\le c\,\hat g_0$ for some $c>0$ and all $t\ge0$.
Hence, $g_\infty$ is positive definite.\qed

\section{The Extrinsic Geometric Flow}
\label{sec:1-def}

We study two types of evolution of Riemannian metrics on $(M,\calf)$, depending on functions $f_m\ (m<n)$, at least one of them is not identically zero:
(a)~$f_m\in C^{2}(M\times\RR)$, (b)~$f_m= f_m(\overrightarrow{\!\tau})\in C^{2}(\RR^{n})$.

\begin{definition}\label{D-EGF}\rm
Given functions $f_m$ of type either (a) or (b), a family $g_t,\ t\in[0,\eps)$, of adapted Riemannian metrics on $(M,\calf, N)$
satisfying along $\calf$ the PDE (\ref{eq1}), see also (\ref{eq1-Lie}), will be called an \textit{Extrinsic Geometric Flow} (EGF).
The~corresponding \textit{normalized EGF} is defined along $\calf$ by
\begin{equation}\label{eq1-n}
 \dt g_t =\nabla^t_N h(b_t)-\frac1n\,r(t)\,\hat g_t,
 \quad{\rm where}\quad r(t)=\int_M N(\tr h(b))\,d\vol_t/\vol(M,g_t).
\end{equation}
\end{definition}

\begin{lemma}\label{P-localA}
Let $g_t$ be the solution to EGF (\ref{eq1}).
Then the Weingarten operator $A$ of $\calf$ w.\,r.\,t. $g_t$ satisfies
\begin{equation}\label{E-Ak2-B}
 \dt A = \sum\nolimits_{m=0}^{n-1}
 \Big(\frac12\,f_m\,[A, \nabla^t_N (A^m)]  -N(N(f_m)) A^m -N(f_m)\nabla^t_N (A^m) -f_m\nabla^t_N\nabla^t_N(A^m)\Big).
\end{equation}
\end{lemma}

\noindent\textbf{Proof}.
Substituting $S^\sharp=\nabla^t_N h(A)$ into (\ref{E-btAt}), see (\ref{eq1}), we obtain (\ref{E-Ak2-B}).\qed

Let $g_t$ be a family of Riemannian metrics of finite volume on $(M,\calf)$.
Metrics $\tilde g_t = (\phi_t\hat g_{t})\oplus g_t^\perp$ with $\phi_t=\vol(M,g_t)^{-2/n}$ have unit volume:
$\int_M d\,\widetilde{\vol}_t=1$.
 The next proposition shows that \textit{unnormalized and normalized EGFs differ only by rescaling along~$\calf$}.

\begin{proposition}\label{P-normEG}
Let $(M,\calf)$ be a foliation, and $g_t$ a solution (of finite volume) to EGF (\ref{eq1}) with
$h(b)=\sum\nolimits_{m=0}^{n-1} f_m(\overrightarrow{\!\tau})\,{b}_m$. Then the metrics
\[
 \tilde g_{\,t}=(\phi_t\,\hat g_t)\oplus g^\perp,\quad{\rm where}\quad
 \phi_t=\vol(M,g_t)^{-2/n},
\]
evolve according to the normalized EGF
\begin{equation}\label{eq1-nB}
 \partial_{\,t}\,\tilde g_t = \nabla^t_N\,h(\tilde b_t) -({\rho_t}/n)\,\hat{\tilde g}_t,\quad
 \mbox{ where }\ \rho_t=\int_M N(\tr h(b))\,d\vol_t\,/\,{\vol(M,g_t)}.
\end{equation}
\end{proposition}

\noindent\textbf{Proof}.
By Remark~\ref{L-btAt}, $\tilde\tau_j=\tau_j$ and $h(\tilde A)=h(A)$ (for metrics $\tilde g_t$ and $g_t$, respectively).
Hence,
\[
 \tr h(\tilde A)=\tr h(A)=\sum\nolimits_{m=0}^{n-1} f_m(\overrightarrow{\!\tau})\,\tau_m.
\]
From (\ref{E-dtdvol}) with $S=\nabla^t_N h(b)$ we get the derivative of the volume function
\[
 \frac{d}{dt}\vol(M,g_t)=
 \frac{d}{dt} \int_M d\vol_t = \frac12 \int_M  N(\tr h(A))\,d\vol_t.
\]
Thus $\phi_t=\vol(M,g_t)^{-2/n}$ is a smooth function.
By Remark~\ref{L-btAt}, we have $h(\tilde b)=\phi_t\cdot h(b)$.
By the above and $\hat{\tilde g}_t=\phi_t\,\hat{g}_t$, we have
\[
 \partial_{\,t}\,\tilde g_t
 =\phi_t\dt g_t + \phi\,'_t\,\hat g_t
 = h(\tilde b) + \phi\,'_t/\phi_t\,\hat{\tilde g}_t.
\]
Using (\ref{E-dtdvol}) and $d\,\widetilde{\vol}_t=\phi^{n/2}_t\,d\vol_t$, we obtain
\[
 \dt\widetilde{\vol}_t=\dt(\phi^{\frac n2}_t \vol_t)
 =\frac n2\,\phi_t^{\frac n2-1}\phi'_t\vol_t
  +\frac12\,\phi^{\frac n2}_t N(\tr h(A))\vol_t
 =\frac12\Big(n\,\frac{\phi'_t}{\phi_t}+N(\tr h(A))\Big)\widetilde{\vol}_t.
\]
Let $\rho_t$ be the average of $N(\tr h(A))$, see (\ref{eq1-nB}). From the above we get
\begin{eqnarray*}
 0\eq 2\,\frac{d}{dt}\int_M d\,\widetilde{\vol}_t =\int_M\Big(n\,\frac{\phi'_t}{\phi_t}
 +N(\tr h(A))\Big)d\,\widetilde{\vol}_t
 = n\,\frac{\phi\,'_t}{\phi_t}+\rho_t.
\end{eqnarray*}
This shows that
 $\rho_t/n = -{\phi\,'_t}/{\phi_t}$.
Hence, $\tilde g$ evolves according to (\ref{eq1-nB}).\qed

Obviously, EGFs preserve Riemannian foliations.
Totally geodesic foliations are the fixed points of (\ref{eq1}).
We shall show that EGFs preserve the umbilicity of~$\calf$. Define the function of one variable~by
\begin{equation*}
 \psi(\lambda)=-\sum\nolimits_{m=0}^{n-1} f_m(n\lambda,n\lambda^2,\ldots,n\lambda^n)\lambda^m.
\end{equation*}
\begin{proposition}\label{P-08b}
Given $g_0$ on $(M,\calf)$, let $g_t\ (0\le t<\eps)$ be a unique smooth solution to the EGF (\ref{eq1}) with $f_m$ of type (b).
If $\calf$ is umbilical for $g_0$ and $\psi'(\lambda)>0$ then $\calf$ is umbilical for any $g_t$.
\end{proposition}

\noindent\textbf{Proof}.
Since $\calf$ is $g_0$-umbilical, we have $A_0=\lambda_0\,\widehat\id$ for some function $\lambda_0: M\to\RR$.
Hence, $h(A_0)=-\psi(\lambda_0)\,\hat\id$.
 Assume that $\calf$ is umbilical for all $g_t$, then
\[
 A_t=\lambda_t\id,\qquad
 h(A_t)=-\psi(\lambda_t)\,\hat\id,
\]
where $\lambda_t: M\to\RR$ are smooth functions.
In this case, by Lemma~\ref{L-btAt2} with $s=-\psi(\lambda_t)$, we have
\begin{equation}\label{E-EGF-tlambda}
 \dt\lambda_t=\frac12\,N(N(\psi(\lambda_t))) =\frac12\,N(\psi'(\lambda_t) N(\lambda_t)).
\end{equation}
Since $\psi'(\lambda)>0$, the PDE (\ref{E-EGF-tlambda}) is parabolic and admits a unique local solution, see Appendix.
Now, (\ref{eq1}) reads~as
\begin{equation}\label{E-surf1}
 \dt{\tilde g}_t =-N(\psi(\lambda_t))\,\hat{\tilde g}_t,
\end{equation}
from which we find
$\hat{\tilde g}_t =\exp(-\int_{s=0}^t N(\psi(\lambda_s))\,ds)\hat g_0$ along $\calf$.
The Weingarten operator of certain metrics $\tilde g_t$ is conformal: $A_t=\lambda_t\id$.
 Thus, ${\tilde g}_t$ is the solution of~(\ref{eq1}).
Our assumptions ensure that a solution is unique, hence ${\tilde g}_t=g_t$ for all~$t$.\qed

\section{Searching for power sums}
\label{subsec:2-3-3}

Although EGF (\ref{eq1}) consists of (second order) non-linear PDEs,
the corresponding power sums $\tau_i\ (i>0)$ satisfy an infinite quasilinear system.

\begin{proposition}\label{L-Tlocalsol2}
Power sums $\{\tau_i\}_{i\in\NN}$ of EGF (\ref{eq1}) with $f_m$ of type (a) satisfy the infinite system
\begin{equation}\label{E-tauAk-tilde}
 \dt\tau_i =-\frac{i}2\sum\nolimits_{m=1}^{n-1}\frac{m f_m}{i+m-1}\,N(N(\tau_{i+m-1})) + a_i,
\end{equation}
where
 $\,a_i=-\frac{i}2\sum\limits_{m=1}^{n-1}[ N(N(f_m))\tau_{i+m-1} -f_m\tr(\nabla^t_N(A^{i-1})\nabla^t_N(A^{m}))]$.
The $n$-truncated system (\ref{E-tauAk-tilde})~is
\begin{equation}\label{E-tauBnm-trunc-a}
 \dt\overrightarrow{\!\tau} =-
 \sum\nolimits_{m=1}^{n-1}\frac m2\,f_m\,(B_{n,1})^{m-1}\,N(N(\overrightarrow{\!\tau})) +a,
\end{equation}
where
$a=a(\overrightarrow{\!\tau},N(\overrightarrow{\!\tau}))$ is a smooth vector-function.
\end{proposition}

\noindent\textbf{Proof}.
Since (\ref{eq1}), we substitute $S=\nabla^t_N h(b)$ into (\ref{E-At-3})$_1$ and obtain
\begin{equation}\label{E-tau-i-a}
  \dt\tau_i =-\frac i2\tr\big(A^{i-1}\nabla^t_N\nabla^t_N\big(\sum\nolimits_{m=0}^{n-1} f_m\,A^m\big)\big),\quad i>0.
\end{equation}
The desired system (\ref{E-tauAk-tilde}) follows from the above and the identity
\begin{equation}\label{E-trAnablaB}
\frac{i+m-1}{m}\,\tr(A^{i-1}\nabla^t_N A^{m})=\tr(\nabla^t_N A^{i+m-1})= N(\tau_{i+m-1}).
\end{equation}
By Proposition~\ref{C-BBnm} (see Appendix), the system (\ref{E-tauAk-tilde}) is equivalent to (\ref{E-tauBnm-trunc-a}).\qed

\begin{proposition}\label{P-tauAk-a}
Power sums $\{\tau_i\}_{i\in\NN}$ of EGF (\ref{eq1}) with $f_m$ of type (b) satisfy the infinite system
\begin{equation}\label{E-tauAk}
 \dt\tau_i= -\frac{i}2\Big[\sum\nolimits_{m=1}^{n-1}\frac{m f_m}{i+m-1}\,N(N(\tau_{i+m-1}))
 +\sum\nolimits_{m=0}^{n-1}\tau_{i+m-1}\sum\nolimits_{s=1}^n f_{m,s}\,N(N(\tau_s))\Big] +a_i,
\end{equation}
where  $N(f_m)=\sum\nolimits_s f_{m,s} N(\tau_s)$ and
\[
  a_i=-\frac{i}2\!\sum\limits_{m=1}^{n-1}\!\!\Big[\frac{2\,m}{i{+}m{-}1}\,N(f_m)N(\tau_{i+m-1})
 -f_m\!\tr(\nabla^t_N(A^{i-1})\nabla_N(A^{m})) +\tau_{i+m-1}\!\!\sum\limits_{a,b=1}^{n}\!f_{m,ab} N(\tau_a)N(\tau_b)
 \Big].
\]
 The $n$-truncated system (\ref{E-tauAk}) is
\begin{equation}\label{E-tauBnm-trunc}
 \dt\overrightarrow{\!\tau} = -(\widetilde{B}_n+\widetilde{A})\,N(N(\overrightarrow{\!\tau})) + a,
\end{equation}
where $a=a(\overrightarrow{\!\tau},N(\overrightarrow{\!\tau}))$ is a smooth vector-function,
the matrices $\widetilde{B}_n$ and $\tilde A$ are defined in (\ref{E-BBB-main}).
\end{proposition}

\noindent\textbf{Proof}. It is similar to proof of Proposition~\ref{L-Tlocalsol2}.
Substituting $S=\nabla_N h(b)$ into (\ref{E-At-3})$_1$, see (\ref{eq1}), we obtain (\ref{E-tau-i-a}).
The desired system (\ref{E-tauAk}) follows from the above and (\ref{E-trAnablaB}).
By Proposition~\ref{C-BBnm}, the system (\ref{E-tauAk}) is equivalent to (\ref{E-tauBnm-trunc}).\qed

\vskip1.5mm\noindent
\textbf{Problem}. Express $a_i\ (i>1)$ in Propositions~\ref{L-Tlocalsol2} and \ref{P-tauAk-a} through the scalar invariants of $A$ and $\nabla^t_{\!N} A$.

\begin{example}\label{Ex-NT}\rm
Among EGFs (\ref{eq1}) with $f_m$ of type (b), the \textit{Newton transformation flow}, see (\ref{eq1-NT1}), corresponds to
$f_m=(-1)^m\sigma_{r-m}$. \textit{Newton transformations} $T_r(A)$ are defined as
\begin{eqnarray}\label{E-Tr1}
 T_r(A)=\sum\nolimits_{i=0}^r(-1)^i \sigma_{r-i}A^i=\sigma_r\,\widehat\id-\sigma_{r-1}\,A+\ldots+(-1)^r A^r\qquad (0\le r\le n),
\end{eqnarray}
see \cite{rw-m}. Since $\tr T_r(A)=(n-r)\,\sigma_r$, the corresponding to (\ref{eq1-NT1}) normalized flow is
\begin{equation}\label{eqNT-n}
 \dt g_t = \nabla^t_N\,T_r(b_t)-\frac1n\,r(t)\,\hat g_t,\quad{\rm where}\quad
 r(t)=(n-r)\int_M N(\sigma_r)\,d\vol_t/\vol(M,g_t).
\end{equation}
For EGF (\ref{eq1-NT1}), denoting $T_j=T_j(A)$, by Lemma~\ref{L-btAt2}, we have
\begin{eqnarray*}
 \dt\sigma_i =-\frac12\big[N(\tr(T_{i-1}\cdot\nabla^t_N\,T_r)) -\tr((\nabla^t_N\,T_{i-1})(\nabla^t_N\,T_r))\big]
 =\frac12\tr((\nabla^t_N\,T_{i-1})(\nabla^t_N\,T_r))\\
 -\frac12\sum\nolimits_{j=0}^r(-1)^j\Big[
 N\big(\sigma_{r-j}\tr(T_{i-1}\cdot\nabla^t_N\,(A^j)) +N(\sigma_{r-j})\tr(T_{i-1}\cdot A^j)\big)\Big].
\end{eqnarray*}
One may express the above (by induction) through $\sigma$'s and their $N$-derivative.
 \end{example}

\section{Local existence and uniqueness of metrics}
\label{sec:proofs}

We shall prove our results concerning the existen\-ce/uniqueness of EGFs (\ref{eq1}) and their corollaries.
For brevity we shall omit the index $t$ for time-depen\-dent tensors $A$, $a,\,b_j$ and functions $\tau_i, \sigma_i$.

Generally, EGFs with $f_m$ of type (a) are solvable under additional conditions (e.g., for $A$ and the auxiliary functions $f_m$).
The following theorem concerns EGFs with $f_m$ of type (a) and is essential in the proof of Theorem~\ref{T-02}
about EGFs  with $f_m$ of type (b).

The~coordinate system described in the following lemma, see \cite[Section 5.1]{cc1},
is here called a \textit{biregular foliated chart}.

\begin{lemma}\label{L-adapted-coords}
Let $M$ be a differentiable manifold with a codimension-1 foliation $\calf$ and a vector field $N$ transversal to $\calf$.
Then for any $q\in M$ there exists a coordinate system $(x_0,x_1,\ldots x_n)$ on a neighborhood $U_q\subset M$ (centered at $q$) such that the leaves on $U_q$ are given by $\{x_0=c\}$ (hence the coordinate vector fields $\partial_i=\partial_{x_i},\ i\ge1$, are tangent
to leaves), and $N$ is directed along $\partial_0=\partial_{x_0}$ (one may assume $N=\partial_0$ at $q$).
\end{lemma}

If $(M, \calf , g)$ is a foliated Riemannian manifold and $N$ is the unit normal to $\calf$, then
in biregular foliated coordinates $(x_0, x_1,\ldots, x_n)$, one may present the metric $g$ in the form
\begin{equation}\label{E-G0ij-g1}
 g=g_{00}\,dx_0^2+\sum\nolimits_{i,j>0}g_{ij}\,dx_i dx_j,
\end{equation}
where one may assume $g_{00}=0$ along the axis $(x_0, 0,\ldots, 0)$.
Let $g^{ij}$ be the entries of the matrix inverse to $(g_{ij})$ and $g_{ij,k}$ the derivative of $g_{ij}$ in the direction of $\partial_k$.
We have $b_{ij}=g(\nabla_{\partial_i}\partial_j, N)$.

\begin{lemma}[see \cite{rw-m}]\label{lem-biregG}
For a metric (\ref{E-G0ij-g1}) in biregular foliated coordinates (of a codimension-one foliation $\calf$) on $(M,g)$, 
one has $ N =\partial_0/\sqrt{g_{00}}$ (the unit normal) and
\begin{eqnarray*}
 \Gamma^j_{i0} \eq (1/2)\sum\nolimits_{\/s} g_{is,0}g^{sj},\quad
 \Gamma^{i}_{00} = -(1/2)\sum\nolimits_{\/s} g_{00,s}g^{si},\quad
 \Gamma^{0}_{ij} = -g_{ij,0}/(2\,g_{00}),\\
 b_{ij}\eq \Gamma^0_{ij}\sqrt{g_{00}}=-\frac12\,g_{ij,0}/\sqrt{g_{00}}
 \quad\mbox{(the second fundamental form)},\\
 A^j_{i}\eq -\Gamma^j_{i0}/\sqrt{g_{00}}=\frac{-1}{2\,\sqrt{g_{00}}}\sum\nolimits_{\/s} g_{is,0}\,g^{sj}
 \quad\mbox{(the Weingarten operator)},\\
 (b_m)_{ij}\eq A^{s_1}_{s_2} \dots A^{s_{m-1}}_{s_{m}} A^{s_{m}}_i g_{js_1}
 \quad\mbox{(the $m$-th ``power" of $\,b_{ij}$)},\\
 \tau_i\eq \Big(\frac{-1}{2\sqrt{g_{00}}}\Big)^i\sum\nolimits_{\,\{r_a\},\, \{s_b\}}
 g_{r_1 s_2,0}\ldots g_{r_{i-1} s_i,0}\,g^{s_1 r_1}\ldots g^{s_i r_i}.
\end{eqnarray*}
In particular, $ \tau_1=\frac{-1}{2\sqrt{g_{00}}}\sum\limits_{\,r,\,s} g_{r s,0}\,g^{r s}$.
\end{lemma}

\noindent\textbf{Proof of Theorem~\ref{T-01}}.
Given $q\in M$, by Lemma~\ref{L-adapted-coords}, there exist biregular foliated coordinates $(x_0,x_1,\ldots x_n)$ on $U_q\subset M$ (with center  at $q$) and the metric has the form (\ref{E-G0ij-g1}).
Then the unit normal to $\calf$ on $U_q$ is $N=g_{00}^{-1/2}\partial_0$.
 Denote $F_{ij}:= h(b)(t, x, g_{ab}, \psi_{ab})_{ij}$ and $\psi_{ab}=\partial_0 g_{ab}$.
The EGF system (\ref{eq1}) with $f_m$ of type~(a) has the following form along $N$-curve 
$\gamma: x\mapsto\gamma(x)$ (on $U_q$):
\begin{eqnarray}\label{E-local2-1}
\nonumber
 \dt g_{ij}
 \eq g_{00}^{-1/2}\big(\partial_0 F_{ij} -\Gamma^k_{0i}F_{kj} -\Gamma^k_{0j}F_{ik}\big)\\
 \eq g_{00}^{-1/2}\big(\partial_0 F_{ij}
 -\frac12\sum\nolimits_{s,k}(\psi_{is}g^{sk}F_{kj}+\psi_{is}g^{sk}F_{ik})\big)\qquad (i, j = 1, \ldots , n),
\end{eqnarray}
where $g_{00}=1$ along the axis $\gamma_0=(x_0, 0,\ldots, 0)$.
Certainly, EGF system (\ref{eq1-Lie}) has a shorter form
\begin{equation}\label{E-local2-Lie}
 \dt g_{ij} = g_{00}^{-1/2}\partial_0 F_{ij}(t, x, g_{ab}, \psi_{ab}).
\end{equation}
For example, if $f_m=0$ for $m\ge2$, then $F_{ij}=f_0 g_{ij}+f_1 b_{ij}\ (i,j>0)$,
and (\ref{E-local2-1}) reads as:
\begin{equation*}
 \dt g_{ij} = -\frac12\,g_{00}^{-1/2}f_1\,\partial^2_{00}\,g_{ij}
 +a_{ij}(t, x, g_{ab}, \psi_{ab}),
\end{equation*}
which is parabolic when $f_1<0$, see also (\ref{E-cond-a}).
Using Theorem~A completes the proof in this case.

Now let $f_m\ne0$ for some $m\ge2$ (e.g., general $f_m$).
Computing the derivative
\[
 \partial_0 F_{ij}=\sum\nolimits_{a\le b}\frac{\partial F_{ij}}{\partial\psi_{ab}}\,\partial^2_{00} g_{ab}
 +\sum\nolimits_{a\le b}\frac{\partial F_{ij}}{\partial g_{ab}}\,\psi_{ab}+\frac{\partial F_{ij}}{\partial x_{0}},
\]
we rewrite (\ref{E-local2-1}) as the quasi-linear system
\begin{equation}\label{E-pF1}
 \dt g_{ij} =g_{00}^{-1/2}\sum\nolimits_{a\le b}\frac{\partial F_{ij}}{\partial\psi_{ab}}\,\partial^2_{00} g_{ab}
 + c_{ij}(t, x, g_{ab}, \psi_{ab}).
\end{equation}
In general, the following square matrix of order $\frac12\,n(n+1)$ is not symmetric:
\[
 d_{\psi}F=\Big\{\frac{\partial F_{ij}}{\partial\psi_{ab}}\Big\}\qquad (i\le j,\ \ a\le b).
\]
We claim that (\ref{E-pF1}) is strong parabolic. If we change the local coordinate system on $M$,
the components $F_{ij}\ (i\le j)$ and $\psi_{ab}\ (a\le b)$ at $q$ will be transformed by the same tensor low.
Hence, $d_{\psi}F$ is a $(1,1)$-tensor on the vector bundle of symmetric $(0,2)$-tensors on $T\calf$.
So,  $d_{\psi}F(q)$ can be seen as the linear endomorphism of the space of symmetric $(0,2)$-tensors on $T_q\calf$.
The order of indices in matrix representation of $d_{\psi}F$ is $[1,1],[1,2],\dots,[1,n],[2,2],[2,3],\dots,[2,n],\dots,[n,n]$.

The strong parabolicity is a pointwise property, so can be considered at any point $q\in M$ in a special biregular foliated chart around $q$,
for example, such that $g_{ij} = \delta_{ij}$ and $A_{ij} = k_i\delta_{ij}$ at $q$ for $t=0$.
(Recall that $(k_i)$ are the principal curvatures of $\calf$).
In this chart, our calculations show that the matrix $d_{\psi}F$ is diagonal, so has real eigenvalues at $q$, and
$$
 \frac{\partial F_{ij}}{\partial\psi_{ab}}=\sum\nolimits_{m\ge1} f_m(q,0)\,\mu_{m;ij}\,\delta^{\{i,j\}}_{\{a,b\}}.
$$
By (\ref{E-cond-a}), the system (\ref{E-pF1}) is strong parabolic: the constant in the condition (\ref{E-PDE-1b}),
see Appendix, is
\[
 c=\inf_{q}\inf_{i,j}\big\{\sum\nolimits_{m=1}^{n-1} f_m(q,0)\,\mu_{m;ij}\big\}>0.
\]
By Theorem~A (in  Appendix), given $q\in M$ there exists a unique solution to (\ref{E-pF1}) which is defined in $U_q$ along the $N$-curve through $q$ for some time interval $[0,\eps_q)$ and satisfies the initial conditions $g_{ij}(0,x)=g_{0, ij}$.
 By Theorem~A, the value $\eps_q$ continuously depends on $q\in M$.
 The claim follows from the above and compactness of $M$.
 \qed

\begin{example}[see \cite{rw-m}]\rm
For $h(b)=b_2$, the matrix $d_{\psi}F(q)$ in an orthonormal frame at any point is
\[
\frac{\partial (b_2)_{ij}}{\partial\psi_{ab}}=
\frac14
\left[\begin{smallmatrix} 2\,\psi_{11}&2\,\psi_{12}&2\,\psi_{13}&0&0&0\\ \noalign{\medskip}\psi_{12}&\psi_{11}+\,\psi_{22}&\psi_{23}&\psi_{{12}}&\psi_{13}&0\\ \noalign{\medskip}\psi_{13}&\psi_{23}&\psi_{11}+\,\psi_{33}&0&\psi_{12}&\psi_{13}\\ \noalign{\medskip}0&2\,\psi_{12}&0&2\,\psi_{22}&2\,\psi_{23}&0\\ \noalign{\medskip}0&\psi_{13}&\psi_{{12}}&\psi_{23}&\psi_{22}+\,\psi_{33}&\psi_{23}\\
\noalign{\medskip}0&0&2\,\psi_{13}&0&2\,\psi_{23}&2\,\psi_{33}
\end{smallmatrix}\right]
\]
with the order of indices $[1,1],[1,2],[1,3],[2,2],[2,3],[3,3]$.
At $q$ for $t=0$ (i.e., $\psi_{ab}=0,\ a\ne b$ and $\psi_{aa}=k_a$) it is diagonal with the elements $\mu_{2;ab}=k_a+k_b$.
 For instance, let $A_0=[\frac12, 1,1, \frac12, 1,\frac12]$ be the diagonal matrix of order six.
 Then the matrix $A_1=A_0\,d_\psi (b_2)$ is symmetric, and the system (\ref{E-pF1}) for $h(b)=b_2$ is ``symmetrizable":
 $A_0\dt g_{ij} = A_1\partial^2_{00} g_{ij} +\mbox{\{low-order terms\}}$.
\end{example}

Remark that for EGFs with general $f_m$ of type (b), the solution procedure requires additional assumptions
(e.g., for the auxiliary functions $f_m$ and the matrix $\widetilde{B}_n$).

\vskip2mm\noindent\textbf{Proof of Theorem~\ref{T-02}}.
Let $A_0$ and $\overrightarrow{\!\tau}^{\,0}$ be the values
of extended Weingarten operator and power sums of the principal curvatures $k_i$ of (the leaves of) $\calf$
determined on $(M,\calf)$ by a given metric $g_0$.

\vskip1.5mm
 a) \textit{Uniqueness}. Let $g^{(1)}_t$ and $g^{(2)}_t$ be two solutions to (\ref{eq1}) with the same initial metric $g_0$.
 Functions $\overrightarrow{\!\tau}^{\,t,1}$ and $\overrightarrow{\!\tau}^{\,t,2}$, corresponding to $g^{(1)}_t$ and $g^{(2)}_t$, satisfy (\ref{E-tauAk}) and have the same initial value $\overrightarrow{\!\tau}^{\,0}$.
 By Theorem~A (see Appendix),
 a solution of the parabolic system (\ref{E-tauAk}) is unique, and we have
 $\overrightarrow{\!\tau}^{\,t,1}=\overrightarrow{\!\tau}^{\,t,2}=\overrightarrow{\!\tau}^{\,t}$ on some positive time interval $[0,\eps_1)$.
 Hence, both $g^{(1)}_t$ and $g^{(2)}_t$ satisfy EGF equation (\ref{eq1}) with $f_m$ of type (a) with known coefficients
 $\tilde f_j(p,t):=f_j(\overrightarrow{\!\tau}^{\,t}(p))$.
 By Theorem~\ref{T-01}, $g^{(1)}_t=g^{(2)}_t$ on some positive time interval $[0,\eps_2)$.

\vskip1.5mm
 b) \textit{Existence}.
 By Proposition~\ref{P-tauAk-a} and Theorem~A, (\ref{E-tauAk}) admits a unique solution $\overrightarrow{\!\tau}^{\,t}$ on a~positive time interval $[0,\eps_1)$.
 By Theorem~\ref{T-01}, (\ref{eq1}) with $f_m$ of type (a) with known functions $\tilde f_j(\cdot,t):=f_j(\overrightarrow{\!\tau}^{\,t},t)$ has a unique solution $g^*_t$ ($g^*_0=g_0$) for $0\le t<\eps^*$.
 The Weingarten operator $A^*_t\ (A^*_0=A_0)$ of $(M,\calf,g^*_t)$ satisfies (\ref{E-Ak2-B}),
 hence the power sums of its eigenvalues, $\overrightarrow{\!\tau}^{\,t,*}\ (\overrightarrow{\!\tau}^{\,0,*}=\overrightarrow{\!\tau}^{\,0})$, satisfy (\ref{E-tauAk-tilde}) with $f_m$ of type (a) with the same coefficient functions $\tilde{f}_j$.
 By Proposition~\ref{L-Tlocalsol2} and Theorem~A, the solution of the problem is unique,
 hence $\overrightarrow{\!\tau}^{\,t}=\overrightarrow{\!\tau}^{\,t,*}$,
 i.e., $\overrightarrow{\!\tau}^{\,t}$ are power sums of eigenvalues of $A^*_t$.
 Finally, $g^*_t$ is a solution to (\ref{eq1}) such that
 $\overrightarrow{\!\tau}^{\,t}$ are power sums of the principal curvatures of $\calf$ in this metric.
\qed

\section{Proof of Theorem~\ref{C-nablab} and Propositions~\ref{C-twisted}--\ref{C-tauf}}
\label{sec:appl}

\noindent\textbf{Proof of Theorem~\ref{C-nablab}}.
In biregular foliated coordinates, (\ref{eq1-1}) has the form of parabolic PDEs
\[
 \dt g_{ij}
 =\frac{1}{\sqrt{g_{00}}}\,\partial_0\Big(\frac{1}{\sqrt{g_{00}}}\,\partial_0 g_{ij}\Big)
 =\frac{1}{g_{00}}\,\partial_{00}\,g_{ij} +\frac{1}{2}\,\partial_0\Big(\frac{1}{g_{00}}\Big)\,\partial_0 g_{ij},
\]
see also (\ref{E-local2-Lie}).
One may assume $g_{00}=1$ along the central $N$-curve $\gamma_0(s)$, hence
  $\dt g_{ij} =\partial_{00}\,g_{ij}$ along~$\gamma_0$.
This linear system
has a unique global solution $g_{ij}(t)$ for all $t\ge0$,
which converges to $\bar g=(\bar g_{ij})$ as $t\to\infty$ (along any $N$-curve), see Appendix.
Hence there is a global solution $g_t\ (t\ge0)$ on~$M$.
Since $N$-curves compose a $S^1$-fibration, from $\dt(\det g_{ij}) = \frac12(\tr_{\!g} S)\det g_{ij}$ for $S=-2\,\mathcal{L}_N\,b$, see (\ref{E-dtdvol}), we conclude that $\bar g$ is positive definite and the limit $\bar g=\lim\limits_{t\to\infty} g_t$ is a smooth metric on~$M$.

The corresponding system for $\tau$'s, see (\ref{E-tauBnm-trunc-a}) with $\widetilde{B}_n =\frac12\,\hat\id$, is diagonal parabolic,
\begin{equation*}
 \dt\tau_i =N(N(\tau_{i}))-i\tr(\nabla_N(A^{i-1})\nabla_N(A)).
\end{equation*}
In particular, we have $\dt\tau_1=N(N(\tau_1))$;
from this easily follows that $N(\tau_1)$ satisfies the heat equation $\dt(N(\tau_1))=N(N(N(\tau_1)))$.
Hence, the limit $N(\bar{\tau}_1)=const$ along $N$-curves. Since $\bar\tau_1$ is bounded on a compact manifold $M$,
we obtain $N(\bar{\tau}_1)=0$.
By the integral formula $\int_M [N(\bar\tau_1)-\bar\tau_1^2]\,d\,\overline{\vol}=0$, see \cite{rw-m}, we conclude that $\bar{\tau}_1=0$ on $M$.
\qed

\vskip2mm
\noindent\textbf{Proof of Proposition~\ref{C-twisted}}.
Define $\phi:M\to\RR$ by the equality $f=e^{\,\phi}$. We~have $\tau_1=-N(\phi)$.
 The evolution $\dt g(t)= (s(t)\,g_1(t))\oplus 0$ preserves twisted product structure.
Denoting $\hat g_t=g_1(t)$, we obtain $g(t)=(f(t)^2\hat g_1)\oplus g_2$,
where $g_2$ does not depend on $t$.
Assuming $f(t)=e^{\phi(t)}$, we find $\hat g_t=e^{2\phi(t)}\hat g_0$.
In this case, $\dt g(t)=2 e^{2\phi(t)}\dt\phi(t)\hat g_0$. Hence, $\dt\phi(t)=\frac12\,s(t)$.
 For $s=-\frac2n N(\tau_1)$, EGF equation (\ref{eq1-1}) reads~as
\begin{eqnarray}\label{E-polynin1}
 \dt\phi \eq (1/n)\,N(N(\phi)).
\end{eqnarray}
Replacing $t\to t/n$, one may reduce (\ref{E-polynin1}) to the heat equation.
Hence, see Proposition~\ref{P-normEG}, there is a unique solution $\phi(t,x,y)$ for all $t\ge0$,
satisfying $\lim\limits_{t\to\infty}\phi(t,x,y)=\bar\phi(x)$. Certainly, $\bar f=e^{\bar\phi}$.
\qed

\vskip2mm\noindent\textbf{Proof of Proposition~\ref{C-ftau}}.
By Remark~\ref{R-FPsi}, we have $F_k'(\tau_1)=\frac kn F_{k-1}(\tau_1)$.
By Proposition~\ref{P-conf} with $s=-N(f(\overrightarrow{\!\tau}))$, we find
\begin{equation*}
 N(s)
 =-\sum\nolimits_{k=1}^n N\big(f_{,k}(\overrightarrow{\!\tau})F'_k(\tau_1)N(\tau_1)\big)
 = -\frac kn\sum\nolimits_{k=1}^nN\big(f_{,k}(\overrightarrow{\!\tau})F_{k-1}(\tau_1)N(\tau_1)\big).
\end{equation*}
By this and Lemma~\ref{L-btAt2}, we conclude that $\tau_1$ satisfies the quasi-linear PDE
\begin{equation}\label{E-tau1-prop3}
 \dt\tau_1=N\big(a(\overrightarrow{\!\tau})\,N(\tau_1)\big),\quad {\rm where}\quad
 a=\frac12\sum\nolimits_{k=1}^n k\,f_{,k}(\overrightarrow{\!\tau})F_{k-1}(\tau_1).
\end{equation}
This PDE is parabolic when $a>0$ (see Section~\ref{subsec:parabPDE}) that  proves (i). The proof of (ii) is similar.
\qed

\vskip2mm\noindent\textbf{Proof of Proposition~\ref{C-tauf}}.
By Lemma~\ref{L-btAt2} with $S=s\,\hat g$ and $s=-N(\tau_1-F)$,
we conclude that the functions $\tau_1-F$ and $N(\tau_1-F)$ satisfy the heat equation (see Appendix) along $N$-curves,
\begin{eqnarray*}
  \dt(\tau_1-F)\eq\dt\tau_1= -\frac n2\,N(-\frac2n\,N(\tau_1-F))=N(N(\tau_1-F)),\\
  \dt(N(\tau_1-F))\eq\dt(N(\tau_1))=N(\dt\tau_{1})=N(N(N(\tau_{1}-F))).
\end{eqnarray*}
There is a unique solution, $\tau^t_1$, for all $t\ge0$; moreover, $\tau_1^t\to\bar\tau_1$
uniformly as $t\to\infty$, with the property $N(\bar\tau_1-F)$ is constant along $N$-curves.
Since $M$ is compact, we find that $N(\bar\tau_1-F)=0$, hence $\bar\tau_1-F$ is constant along $N$-curves.
 With known $\tau^t_1$, (\ref{eq1-1f}) also has a unique global solution
 $\hat g_t=\hat g_0\exp(-\!\int_0^t N(\tau^s_1-F)\,ds)$ for all $t\ge0$.
Let us apply the normalized EGF
\begin{equation*}
 \dt\tilde g_t = -\big[\frac2n\,N(\tau_1-F)+\frac2n\,r(t)\big]\hat{\tilde g}_t,\qquad
 r(t) = -2\int_M N(\tau_1-F)\,d\,\widetilde{\vol}_t/\vol(M, \tilde g_t).
\end{equation*}
By Remark~\ref{L-btAt}, the mean curvature is the same for both EGFs: $\tau_1^t=\tilde\tau_1^t$.
Let us differentiate the functional $J_F(t)=\int_M F\,d\,\widetilde{\vol}_{t}$, using
(\ref{E-dtdvol}) with $s_t=-\frac2n\,N(\tau_1-F)-\frac2n\,r(t)$:
\[
 J'_F(t)=-\int_M N(\tau_1-F)\,d\,\widetilde{\vol}_{t} -\frac12\,r(t)\vol(M, \tilde g_t) =0.
\]
By conditions, $J_F(0)=0$. Hence $J_F(t)\equiv0$, i.e., $F$ has zero average on $(M,\,\widetilde{g}_t)$.
 Since $N$-curves are dense in $M$, we have $\bar\tau_1-F=C$  for some $C\in\RR$. Thus, $\tau_1^t-F=C+o(1/t)$.
Integrating over $(M,\,\widetilde{g}_t)$, we obtain
\[
 [C+o(1/t)]\vol(M, \widetilde g_t)=\int_M (\tau_1^t-F)\,d\,\widetilde{\vol}_t=\int_M\tau_1^t\,d\,\widetilde{\vol}_t.
\]
Since $\int_M\tau_1^t\,d\,\widetilde{\vol}_t=0$ (the $\widetilde g_t$-average of $\tau_1^t$ is zero),
see (\ref{E-IF-Reeb}), and $\vol(M, \widetilde g_t)=1$, we have $C=0$.
\qed

\vskip1.5mm\noindent\textbf{Proof of Theorem~\ref{C-tau-sigma-k}}.
By Proposition~\ref{C-ftau}, $\tau_1$ satisfies parabolic PDE (\ref{E-tau1-prop3}).
By Proposition~\ref{P-heat-eq}, there is a global solution $\tau_1^t\ (t\ge0)$ and
we have $\tau_i^t=F_{i}(\tau_1^t)$ for $i>1$.
The metrics $g_t\ (t\ge0)$ are given by the formula $g_t=g_0\,e^{-\int_0^t f(\overrightarrow{\!\tau}^s)\,ds}$.
Moreover, there is a limit $\tau_1^\infty$ (as $t\to\infty$) which is constant along $N$-curves.
By integral formula $\int_M [N(\tau_1)-(\tau_1)^2]\,d\vol=0$, see~\cite{rw-m}, applying for $g_t$, we have $\tau_1^\infty=0$ on $M$. The rate of converging $\tau_i^t\to F_i(0)\ (i>0)$ is exponential, see Proposition~\ref{P-heat-eq}.
By~Proposition~\ref{P-converge}, see also the proof of Theorem~\ref{C-nablab}, there is a limit metric $g_\infty=\lim\limits_{t\to\infty} g_t$.
The~case (ii) is similar.
\qed

\section{EGF on foliated surfaces}
\label{subsec:2-8-3}

 In this section, $(M^2, g_0)$ is a surface equipped with a transversally orientable foliation $\calf$ by curves.
Given $\psi\in C^3(\RR)$, let a 1-parameter family of Riemannian metrics $g_t\ (0\le t<\eps)$ satisfies
the PDE
\begin{equation}\label{E-surf2}
 \dt g_t =-N(\psi(\lambda_t))\,\hat g_t
\end{equation}
see (\ref{eq1-1a}) and (\ref{E-surf1}).
Applying (\ref{E-dtdvol}) with $S=-N(\psi(\lambda_t))\,\hat g_t$, we have
 $\frac{d}{dt}\,(d\vol_t)= -\frac12\,N(\psi(\lambda_t))\,d\vol_t$.
Hence, the volume of a closed surface $(M^2, g_t)$ satisfies the equation
\begin{equation}\label{E-intpsi}
 \frac{d}{dt}\vol(M,g_t)=-\frac12\int_M N(\psi(\lambda_t))\,d\vol_t
 =-\frac n2\int_M \lambda_t\,\psi(\lambda_t)\,d\vol_t.
\end{equation}

\begin{proposition}\label{P-surf}
If $0<C_0(\lambda)\le\psi'(\lambda)\le C_1(\lambda)<\infty$, then (\ref{E-surf2})
on a surface $(M^2,\calf)$ has a unique global smooth solution $g_t\ (t\ge0)$.
\end{proposition}

\noindent\textbf{Proof}.
The $g_t$-geodesic curvature $\lambda_t$ of $\calf$ (w.\,r.\,t. $N$ -- the $g_t$-unit normal to $\calf$)
satisfies the PDE
\begin{equation}\label{E-EGF-tlambda-surf}
 \dt\lambda=\frac12\,N(\psi'(\lambda) N(\lambda)),
\end{equation}
see (\ref{E-EGF-tlambda}).
By Proposition~\ref{P-heat-eq} in Section~\ref{subsec:parabPDE},
(\ref{E-EGF-tlambda-surf}) (with $\lambda_0$ determined by $g_0$) admits a unique smooth solution $\lambda_t$ for all $t\ge0$.
In this case, the EGF (\ref{E-surf2})
has a unique smooth solution $g_t$ on a surface $M^2$ for all $t\ge0$, moreover,
$\hat g_t=\hat g_0\exp(-\int_0^t N(\psi(\lambda_t))\,dt)$.\qed

\begin{example}\label{ex-simple}\rm
Let $\psi=2\,\lambda$, i.e., $\dt g=-N(2\,\lambda_t)\,\hat g_t$, and conditions of Proposition~\ref{P-surf} are satisfied.
Then (\ref{E-EGF-tlambda-surf}) is the heat equation $\dt\lambda_t=N(N(\lambda_t))$ along $N$-curves,
whose solution is known.
The EGF solution is $\hat g_t=\hat g_0\exp\big(-\int_0^t N(\lambda_\xi)\,d\xi\big)\ (t\ge0)$.
For closed manifold $M$, by~(\ref{E-intpsi}), we have $\frac{d}{dt}\vol(M,g_t)=\int_0^t \lambda^2_\xi\,d\xi\le0$,
hence the function $\vol(M,g_t)$ is decreasing.
\end{example}

\begin{proposition}\label{L-nablaNNt}
The Gaussian curvature of the EGF (\ref{E-surf2}) on $(M^2,\calf)$ is given by formula
\begin{equation}
\label{E-gaussM2}
 K_t = \Div\Big(\exp\big(\int_0^t N(\psi(\lambda_\xi))\,d\xi\big) \nabla^0_NN\Big) +N(\lambda_t) -\lambda_t^2.
\end{equation}
\end{proposition}

\noindent\textbf{Proof}.
Notice that $\tau_j=n\,\lambda^j$.
By Lemma~\ref{L-nablaNN} with $S=-N(\psi(\lambda_t))\,\hat g_t$ (and $n=1$), we have
 $\dt(\nabla^t_NN)=N(\psi(\lambda_t))\nabla^t_NN$.
Integrating the above ODE~yields
 \begin{equation}
\label{E-nablaNNt}
  \nabla^t_NN = \exp\Big(\int_0^t N(\psi(\lambda_\xi))\,d\xi\Big) \nabla^0_NN.
 \end{equation}
In our  case here, $K=\Ric(N,N)$, therefore (\ref{E-gaussM2}) is a consequence of (\ref{E-nablaNNt}) and
the identity for codimension-one foliations (see \cite[Proposition~1.4]{rw-m} with $r=0$),
\[
 {\rm div}(\nabla_N N +\tau_1 N) =\Ric(N,N) + \tau_2 -(\tau_1)^2.
\]
Indeed, one should apply the identities $\Div N =-\tau_1$ and $\Div(\tau_1 N) =N(\tau_1)-(\tau_1)^2$.\qed

\begin{example}\label{ex-Reeb1}\rm
Let a function $f\in C^2(-1,1)$ has vertical asymptotes $x=\pm1$ and exactly one strong minimum at $x=0$.
For example, $f=\frac1{10}\big(e^{x^2/(1-x^2)}-1\big)$.
Consider a \textit{Reeb foliation} $\calf$ in the strip $\Pi=[-1,1]\times\RR$ (equipped with the flat metric $g_0$)
whose leaves are
\[
 L_{\pm}=\{x=\pm1\},\quad
 L_s(x)= \{(x,\ f(x)+s),\ |x|<1\},\quad {\rm where}\quad s\in\RR.
\]
Identifying the points $(-1,y)\simeq (1,y)$, we obtain a foliated cylindrical surface $S^1\times\RR$.
Next, identifying the points $(1,y)\simeq (1,y+n)$ for $n\in\ZZ$, we obtain a foliated flat torus $\Pi'=S^1\times S^1$.

Denote $\alpha(x)$ the angle between the leaves $L_s$ and the $x$-axis at the intersection points. That is, $f$ and $\alpha$ are related~by
\[
 f'(x)=\tan\alpha(x)\quad\Leftrightarrow\quad
 \cos\alpha=[1+(f')^{2}]^{-1/2},\quad \sin\alpha=f'\,[1+(f')^{2}]^{-1/2}.
\]
The geodesic curvature of $L_s$ is
\[
 \lambda_0(x)={f''(x)}[1+(f'(x))^{2}]^{-3/2}=\alpha'(x)\cdot|\cos\alpha(x)|,\quad |x|<1.
\]
The tangent and unit normal to $\calf$ vector fields (on the strip) are
 $X=[\cos\alpha(x),\,\sin\alpha(x)],\ N=[-\sin\alpha(x),\,\cos\alpha(x)]$.
The unit normal $N$ at the origin is directed along $y$-axis.
 We use $X$ and $N$ to represent the standard frame
 $e_1=\cos\alpha(x)X-\sin\alpha(x)N,\ e_2=\sin\alpha(x)X+\cos\alpha(x)N$ in $\Pi$.

Taking $\alpha(x)=\frac\pi2 x$, we obtain
the Reeb foliation with $\lambda_0(x)=\frac\pi2\cos(\frac\pi2\,x)>0$.

Let the Riemannian metrics $g_t$ for $t\in[0,\eps)$ satisfy (\ref{E-surf2}) on $\Pi'$.
 Hence $\hat g_t=\hat g_0\,e^{-U_t}$, where
\[
 U_t(x)=\int_0^t N(\psi(\lambda_\xi(x)))\,d\xi
 =-\sin\alpha(x)\int_0^t\dx\psi(\lambda_\xi(x))\,d\xi.
\]
 So,
 $g_t(X,X)=e^{-U_t},\ g_t(X, N)=0$, and $g_t(N, N)=1$.
The scalar products
$g_{11}(t)=g_t(e_1,e_1)$, $g_{12}(t)=g_t(e_1,e_2)$ and $g_{22}(t)=g_t(e_2,e_2)$ are
\begin{equation*}
 g_{11}(t)=\sin^2\alpha+\cos^2\alpha\,e^{-U_t},\quad
 g_{12}(t)=\sin\alpha\cos\alpha(e^{-U_t}-1),\quad
 g_{22}(t)=\cos^2\alpha+\sin^2\alpha\,e^{-U_t}.
\end{equation*}
In particular, $|g_t|=g_{11}(t)g_{22}(t)-g_{12}(t)^2=e^{-U_t}$.
The $N$-curves $\gamma(s)=(x(s),y(s))$ satisfy the system of ODEs
  $\frac{d}{ds} x=-\sin\alpha(x)$, $\frac{d}{ds} y=\cos\alpha(x)$.
From the first of ODEs we deduce the implicit formula
 $s=-\int_x^{\phi_s(x)}\frac{d\xi}{\sin\alpha(\xi)}$
for local diffeomorphisms $\phi_s\ (s\ge0)$  of the interval $(-1,\,1)$.
(Indeed, moving along $N$-curve by the length $s$ from a point $(x,y)$ to a point $(x',y')$, we find that $x'=\phi_s(x)$).

\vskip1.5mm
\underline{Let $\psi(\lambda)=2\,\lambda$}, see Example~\ref{ex-simple}.
By Proposition~\ref{P-surf}, the EGF $\dt g_t = -2\,N(\lambda_t)\,\hat g_t$
admits on $\Pi'$ a unique global solution $g_t\ (t\ge0)$,
and (by Example~\ref{ex-simple}) the function $\vol(\Pi',g_t)$ is decreasing.
 Hence $\hat g_t=\hat g_0\,e^{-U_t}$, where $U_t=-\sin\alpha(x)\int_0^t \dx\lambda_\xi(x)\,d\xi$,
and $\lambda$ evolves along $N$-curves with arc-length parameter $s\in\RR$
by the heat equation $\dt\lambda=\partial^2_{ss}\lambda$.
In terms of $x$, we have $\dt\lambda=\sin\alpha(x)^2\partial^2_{xx}\lambda$.
This heat conduction PDE in a material with variable ``thermal conductivity" $\sin\alpha(x)^2$ on a circle
is parabolic everywhere except one (singular) point $x=0$.
 Using Example~\ref{Ex-parabolic1}\,(i), we write solution along $N$-curves as
\[
 \lambda_t(\phi_0(s))=\int_\RR \lambda_0(\phi_0(\xi)) G(t,s,\xi)\,d\,\xi,\quad
 \mbox{where}\quad G(t,s,\xi)=\frac 1{(4\pi t)^{1/2}}\, e^{-(s-\xi)^2/(4t)}.
\]
Since $\lambda_0(\phi_0(s))\to0$ quickly as $|s|\to\infty$, we have uniform convergence $\lambda_t(x)\to0$ as $t\to\infty$ for any $x\in(-1,1)$.
Along the leaves $L_{\pm\infty}$ the metric does not change, hence $\lambda_t(\pm1)=0$ for all $t\ge0$.

It is known that a curve on a surface with $\lambda=0$ is a geodesic.
From the above we conclude that $g_t$ do not converge as $t\to\infty$ to a Riemannian metric on a torus $\Pi'$,
see Remark~\ref{R-replace}
(otherwise $\calf$ becomes a geodesic foliation, that contradicts to the fact that Reeb foliation is not geodesible).
Since $g_t$ do not depend on the $y$-coordinate, the Gaussian curvature is given by (see \cite{rw-m})
\[
 K_t=-\frac{1}{2\sqrt{|g_t|}}\,\dx\Big(\frac{\dx g_{22}}{\sqrt{|g_t|}}\Big)
 =-\frac12\,e^{U_t}\big(\ddx g_{22}+\frac12\,(\dx U_t)(\dx g_{22})\big).
\]
Substituting
$g_{22}(t)=\cos^2\alpha+\sin^2\alpha\cdot e^{-U_t}$ into the above formula yields
\begin{eqnarray}\label{E-Kiii}
\nonumber
 2\,e^{-U_t} K_t \eq -\big(2\cos(2\alpha)\,(\alpha')^2+\sin(2\alpha)\,\alpha''\big)(e^{-U_t}-1) +\sin^2\alpha\,e^{-U_t}\ddx U_t\\
 \plus\frac32\sin(2\alpha)\,\alpha' e^{-U_t}\dx U_t -\frac12\sin(2\alpha)\,\alpha'\dx U_t -\frac12\sin^2\alpha\,e^{-U_t}(\dx U_t)^2.
\end{eqnarray}
Because $\alpha(0)=0$ and $U_t(0)=0$, by (\ref{E-Kiii}) one has $K_t(0)\equiv0$.
For $\alpha=\frac\pi2 x$ and small $x$ we have
\[
 e^{-U_t} K_t = \frac38\,\pi^3 V_t(0)\,x +o(x),
\]
where $V_t(x)=\int_0^t \dx\lambda_\xi(x)\,d\xi$.
Thus $K_t(x)$ (for $t>0$) changes its sign when we cross the line $x=0$.
\end{example}

\section{Appendix: PDEs}
\label{subsec:2-3-2}

In Section~\ref{subsec:parabPDE}, follow \cite{ta} we discuss parabolic PDEs and prove Proposition~\ref{P-heat-eq} about the quasi-linear heat equation.
In Section~\ref{subsec:GCmatrix}, using the generalized companion matrix, defined in \cite{rw-m}, we prove Proposition~\ref{C-BBnm} about infinite system of parabolic PDEs.

\subsection{Parabolic PDEs in one space variable}
\label{subsec:parabPDE}

Now, we recall some facts about parabolic systems of quasi-linear PDEs.

Let $A=(a_{ij}(t,x,\bu))$ be an $n\times n$-matrix, $a=(a_{i}(t,x,\bu,\dx\bu))$ -- an $n$-vector,
and  $t,x\in\RR$.
Consider the \textit{quasilinear} system of PDEs, $n$ equations in $n$ unknown functions $\bu=(u_1,\ldots ,u_n)$
\begin{equation}\label{E-PDE-1}
 \dt\bu = A(t,x,\bu)\,\ddx\bu + a(t,x,\bu,\dx\bu).
\end{equation}
PDEs (\ref{E-PDE-1}) are \textit{homogeneous} if $a\equiv0$. When $A$ depends on $t$ and $x$ only, the system is \textit{semilinear}. When $A$ and $a$ are functions of $t$ and $x$, the system is \textit{linear}.
The \textit{initial value problem} for (\ref{E-PDE-1}) with given smooth data $A,\,a$ and $\bu_0$,
\begin{equation}\label{E-PDE-ic}
\bu(0,x)=\bu_0(x),
\end{equation}
consists in finding smooth function $\bu(t,x)$ satisfying (\ref{E-PDE-1})\,--\,(\ref{E-PDE-ic}).

The \textit{parabolicity condition} for (\ref{E-PDE-1}) says that
there is $c=const>0$ such that
\begin{equation}\label{E-PDE-1b}
 \<A\,v,v\>\ge c\<v,v\>,\quad \forall\,v\in\RR^n.
\end{equation}

\vskip2mm\noindent
\textbf{Theorem A} \cite[Ch.~15, Proposition 8.3]{ta}.
\textit{
Suppose that $A$ of class $C^\infty$ satisfies the parabolicity condition (\ref{E-PDE-1b}).
Then the system (\ref{E-PDE-1})--(\ref{E-PDE-ic})
with  $\bu_0\in H^s(\RR)$ (the Sobolev space) for some $s>1$,
admits a unique solution
 $u\in C([0,T),\ H^s(\RR))\cap C^\infty((0,T]\times\RR)$,
which persists as long as $\|u\|_{C^{r}}$ is bounded, given $r > 0$.}

\vskip3mm
One may apply this and the maximum principle to obtain
the global existence and uniqueness result for a scalar parabolic equation (i.e., $n=1$ and $A_{11}=k(u)>0$) on $S^1$.

The following proposition is standard, for the convenience of the reader, we give its proof.

\begin{proposition}\label{P-heat-eq}
 Let  a function $k(u)\in C^\infty(\RR)$ (the thermal diffusivity) satisfies
 \[
  0<c_1\le k(u) \le c_2<\infty
 \]
 for some real $c_1\le c_2$. Then the quasi-linear heat equation on a unit circle $S^1$
\begin{equation}\label{E-PDE-1-c1-div}
  \dt u = \dx(k(u)\,\dx\bu),
  \qquad u(0,x)=u_0(x) \in C^\infty(S^1)
\end{equation}
admits a unique solution $u\in C^\infty([0,\infty)\times S^1)$.
Moreover, there exists $\lim\limits_{t\to\infty}u(t,x) =u_\infty\in\RR$, and for some real $\alpha,K>0$ 
the following inequalities are satisfied:
\begin{equation}\label{E-PDE-A}
 \|u(t,\cdot) -u_\infty\|_{S^1}\le  K\,\frac{c_2}{c_1}\,e^{-\alpha\,t}\|u_0 -u_\infty\|_{S^1}.
\end{equation}
\end{proposition}

\noindent\textbf{Proof}.
By  Proposition 9.11 in \cite[Ch.~15]{ta}, there is a unique solution $u\in C^\infty([0,\infty)\times S^1)$.
By the maximum principle, $\|u(t,\cdot)\|_{S^1}\le\|u_0\|_{S^1}$ for all $t>0$.
Define the function $v=\varphi(u)$ of variables $(t,x)$, where $\varphi'(u)=k(u)$.
The monotone function $\varphi$ (of one variable) satisfies inequalities
\begin{equation}\label{E-PDE-varphi}
 c_1\le\varphi'\le c_2.
\end{equation}
The PDE (\ref{E-PDE-1-c1-div})$_1$ reads as $\dt u= \ddx v$, hence (in view of derivation $\dt v=\varphi\,{'}(u)\,\dt u$)
is equivalent to
\[
 \dt v = k(u)\,\ddx v,\qquad v(0,x)=\varphi(u_0(x)).
\]
 Let us compare it with the linear PDE on $S^1$,
\begin{equation}\label{E-PDE-w}
 \dt\widetilde v =\tilde k(t,x)\,\ddx \widetilde v,\qquad \widetilde v(0,x)=\varphi(u_0(x)),
\end{equation}
where $\tilde k(t,x)=k(u(t,x))$ is given. Indeed, $c_1\le \tilde k(t,x)\le c_2$ for all $t>0$ and $x\in S^1$.
From the existence and uniqueness of a solution to (\ref{E-PDE-w}) we conclude that $\widetilde v=v$.
 Denote $u_\infty\in\RR$ the average of $u_0$ over $S^1$, and set $v_\infty=\varphi(u_\infty)$.
The function $w=v(t,x)-v_\infty$ solves the linear PDE on $S^1$
\begin{equation}\label{E-PDE-ww}
 \dt\widetilde w =\tilde k(t,x)\,\ddx\widetilde w,\qquad\widetilde w(0,x)=\varphi(u_0(x))-v_\infty,
\end{equation}

By the theory of linear PDEs, 
(\ref{E-PDE-ww}) possesses a fundamental solution $\widetilde G(t,x,y)$
which can be built by the classical parametrix method, and satisfies the inequalities
$0\le\widetilde G\le K G$ for some real $K>0$. Here $G$ is the fundamental solution of the heat equation $\dt u=\alpha\,\ddx u$ for some constant $\alpha>0$, see survey~\cite{IKO}.
Recall that $G(t,x,y)=\sum_{j\ge0} e^{-\alpha\, j^2\,t}\phi_j(x)\,\phi_j(y)$, where $\phi_j$ denotes the eigenfunction
(of operator $-\alpha\,\ddx$ on $S^1$ with eigenvalue $\lambda_j=\alpha j^2$) satisfying $\int_{S^1}\phi_j(x)\,dx=1$.
Indeed, $\widetilde w(t,x)=\int_{S^1} \widetilde G(t,x,y)\,(\varphi(u_0(y))-v_\infty)\,dy.$
 In particular,
\[
 \|\widetilde v(t,\cdot)-v_\infty\|_{S^1}\le e^{-\alpha\, t}\|\widetilde v(0,\cdot)-v_\infty\|_{S^1}.
\]
Thus, for all $t>0$ and $x\in S^1$ we have a priori estimate
\[
 \|\varphi(u(t,\cdot))-\varphi(u_\infty)\|_{S^1}\le K\,e^{-\alpha\,t}\|\varphi(u_0)-\varphi(u_\infty)\|_{S^1}.
\]
From this, using inequalities
\[
c_1\|u(t,\cdot)-u_\infty\|_{S^1}\le \|\varphi(u(t,\cdot))-\varphi(u_\infty)\|_{S^1},\quad
 e^{-\alpha\,t}\|\varphi(u_0)-\varphi(u_\infty)\|_{S^1}\le c_2\,e^{-\alpha\,t}\|u_0-u_\infty\|_{S^1},
\]
see (\ref{E-PDE-varphi}), we deduce (\ref{E-PDE-A}).\qed

\begin{example}\label{Ex-parabolic1}\rm
(i) Consider the \textit{heat equation} over infinite (or semi-infinite) spatial interval $\RR$,
\begin{equation}\label{E-PDE-1-c2}
  \dt u = \ddx\bu + a(t,x),\quad u(0,x)=u_0(x)\in L^2(\RR).
\end{equation}
Here we assume that $u_0(x)$ is a bounded function.
All the solutions of (\ref{E-PDE-1-c2}),
\begin{equation*}
 u(t,x)=\int_\RR u_0(\xi) G(t,x,y)\,d\,y
 +\int_0^t\int_\RR a(y,s) G(t-s,x,y)\,d\,y\,ds,
\end{equation*}
include the \textit{heat kernel} $G(t,x,y)={(4\pi t)^{-1/2}}\, e^{-(x-y)^2/(4t)}$
(the fundamental solution of homogeneous (\ref{E-PDE-1-c2}) for $u_0(x)=\delta_x(\xi)$ -- the Dirac delta function).
For any function $u_0\in L_2(\RR)$, a unique solution of the homogeneous heat equation,
$u(t,x)=\int_\RR u_0(y)G(t,x,y)\,d\,y$, converges uniformly to a linear function, as $t\to\infty$.
Indeed, if $u_0$ is bounded then the linear function is constant.

(ii) The \textit{rough Laplacian} for tensors is defined by $\Delta=\Div\nabla$.
For a unit circle $S^1$, the eigenvalues of $-\Delta$ are the numbers $\lambda_j=j^2$ with
eigenfunctions $\phi_j(x)=\frac1{\sqrt{2\pi}}e^{-i j x}$, where $j\in\ZZ\setminus\{0\}$.
The~heat equation on $S^1$,
\begin{equation}\label{E-Heat-circle}
  \dt u = \ddx u,\qquad u(0,\cdot)=u_0\in H^2(S^1),
\end{equation}
has a unique solution for functions  $u\in C([0,\infty),\,H^2(S^1))\cap C^1((0,\infty],\,L^2(S^1))$.
By Sobolev embedding theorem, $H^2(S^1)\subset C^1(S^1)$.
The solution of (\ref{E-Heat-circle}) has the property that $u(t,\cdot)\in C^\infty(S^1)$ for all $t>0$.
Moreover, $u(t,\cdot)\to\bar u_0$ as $t\to\infty$, where $\bar u_0=\frac1{2\pi}\int_{S^1} u_0(x)\,dx$
and
\[
 \|u(t,\cdot)-\bar u_0\|\le e^{-t}\|u_0-\bar u_0\|.
\]
Consider the Jacobi \textit{theta function}
\[
 \theta(z,\tau)=\sum\nolimits_{n=1}^\infty  e^{\pi i n^2\tau+2\pi i n z}
 =1+2\sum\nolimits_{n=1}^\infty (e^{\pi i\tau})^{n^2}\cos(2\pi n z),
\]
where $i^2=-1$, $z$ is a complex number and $\tau$ is confined to the upper half-plane.
Taking $z=x\in\RR$ and $\tau=4\pi t\,i$ with real $t>0$, we write
 $\theta(x, 4\pi t\,i)=1+2\sum\nolimits_{n=1}^\infty e^{-4 \pi^2 n^2 t}\cos(2\pi n x)$, which satisfies the heat equation
 (\ref{E-Heat-circle})$_1$.
Since $\lim\limits_{t\to 0}\theta(x,  4\pi t\,i)=\sum\nolimits_{n\in\ZZ}\delta(x-n)$,
the solution of (\ref{E-Heat-circle}) can be specified by convolving the periodic boundary condition at $t=0$ with
$\theta(x, 4\pi t\,i)$.

(iii) Follow \cite{DV}, denote $U(t,x)=\frac{\sin x}{\sqrt{\cos^2 x+e^{2t}}}\ (t\ge0)$ and $k(u)=\frac{1}{1+u^2}$.
We have a 3-parameter family $u(t,x)=U(\omega^2(t+\beta_1), \omega(x+\beta_2))$
of exact solutions to the PDE (\ref{E-PDE-1-c1-div}),
$\dt u=\dx(k(u)\,\dx u))$ on $S^1$.
Set $\omega=1$ and $\beta_1=\beta_2=0$, hence $u(t,x)=U(t,x)$ and $u_0(x)=\frac{\sin x}{\sqrt{\cos^2 x+1}}$.
Indeed, $\lim\limits_{t\to\infty}u(t,x)=u_\infty=0$ for all $x\in S^1$.
Since $0\le u^2\le1$, we conclude that $\frac12\le k(u)\le 1$.
Finally, we have $\|u(t,\cdot)\|_{S^1}\le e^{-t}$ and $\|u_0\|_{S^1}=1$, that is consistent with (\ref{E-PDE-A}).
\end{example}

\subsection{The generalized companion matrix}
\label{subsec:GCmatrix}

Let $P_n=k^n-p_1 k^{n-1}-\ldots-p_{n-1} k-p_n$ be a polynomial over $\RR$ and
$k_1\le k_2\le\ldots\le k_n$ be the roots of $P_n$ for $n>0$.
Hence, $p_i=(-1)^{i-1}\sigma_i$, where $\sigma_i$ are elementary symmetric functions of the roots $k_i$.
The following \textit{generalized companion matrix} (see \cite{rw-m})
plays a key role in this work:
\begin{equation}\label{E-Bn1}
 B_{n,1}=\left(\begin{smallmatrix}
           0 & \frac12 & 0 & \cdots & 0\\
           0 & 0 & \frac23 & \cdots  & 0\\
           \cdots & \cdots & \cdots & \cdots & \cdots\\
           0 & 0 & 0 & \cdots & \frac{n-1}{n} \\
           (-1)^{n-1}\frac{n}{1}\,\sigma_n\quad& (-1)^{n-2}\frac{n}{2}\,\sigma_{n-1}\quad& \ldots\quad& \ -\frac{n}{n-1}\,\sigma_2\ \ & \sigma_1 \\
     \end{smallmatrix}\right).
\end{equation}

\begin{lemma}[see \cite{rw-m}]\label{L-Cn2}
The matrix (\ref{E-Bn1}) has the following properties:

\vskip1.5mm
a) The characteristic polynomial of $B_{n,1}$ is $P_n$.

\vskip1.5mm
b) $v_j=(1,\,2\,k_j,\,3\,k_j^2,\ldots, n\,k_j^{n-1})$ is the eigenvector of $B_{n,1}$ for the eigenvalue $k_j$.

\vskip1.5mm
c) $B_{n,1} V = V D$, where
  $V=\{\frac{n}{i} k_j^{i-1}\}_{1\le i,j\le n}$ is the Vandermonde type matrix, and

\quad $D={\rm diag}(k_1,\ldots, k_{n})$ is a diagonal matrix.
 (If all $k_i$'s are distinct, then $V^{-1}B_{n,1}\,V = D$).
 \end{lemma}

\begin{proposition}[see \cite{rw-m}]\label{L-BBnm}
Let $\tau_i(t,x)\ (i\in\mathbb{N})$ be the power sums
of smooth functions $k_i(t,x)\ (1\le i\le n)$.
Given $m>0$, consider the infinite system of linear PDEs
\begin{equation}\label{E-tauxBnm}
  \dt\tau_i=-\frac{i\,m}{2(i+m-1)}\,\dx\tau_{i+m-1},\qquad i\in\NN.
\end{equation}
Then the $n$-truncated (\ref{E-tauxBnm}),
i.e., $\tau_{n+i}$'s are eliminated using suitable polynomials of $\tau_1,\ldots ,\tau_n$, is
\begin{equation*}
 \dt\overrightarrow{\!\tau} =-
 \frac{m}2\,(B_{n,1})^{m-1}\dx\overrightarrow{\!\tau}.
\end{equation*}
\end{proposition}

Remark that for $m=1$ the system (\ref{E-tauxBnm}) has diagonal form: $\dt\tau_i =-\frac{1}2\,\dx\tau_{i}$,
and for $m=2$ the $n$-truncated system (\ref{E-tauxBnm}) reads as $\dt\tau_i =-B_{n,1}\,\dx\tau_{i}$.

\begin{proposition}\label{C-BBnm}
 Let $\tau_i(t,x)\ (i\in\mathbb{N})$
 be as in Proposition~\ref{L-BBnm}.
Given $m>0$, consider the infinite system of linear PDEs
\begin{equation}\label{E-tauBnm}
  \dt\tau_i=\frac{i\,m}{2(i+m-1)}\,\ddx\tau_{i+m-1},\qquad i\in\NN.
\end{equation}
Then the $n$-truncated system (\ref{E-tauBnm}) has the form
\begin{equation*}
 \dt\overrightarrow{\!\tau} =
 \frac{m}2\,(B_{n,1})^{m-1}\ddx\overrightarrow{\!\tau}+a_m(\overrightarrow{\!\tau},\dx\overrightarrow{\!\tau}).
\end{equation*}
\end{proposition}

\noindent\textbf{Proof}. By Proposition~\ref{L-BBnm}, we have
 $\dx\tau_{n+i}=\sum\nolimits_{j=1}^n\tilde b_{ij}\,\dx\tau_{j}$,
where $\tilde b_{ij}$ are the coefficients of the matrix $B_{n,m-1}$. Derivation of this, leads to
 $\ddx\tau_{n+i}=\sum\nolimits_{j=1}^n\tilde b_{ij}\,\ddx\tau_{j}+\dx\tilde b_{ij}\,\dx\tau_{j}$.\qed

\begin{remark}\rm
For $m=1$ the system (\ref{E-tauBnm}) has diagonal form: $\dt\tau_i =\frac{1}2\,\ddx\tau_{i}$,
and for $m=2$ the $n$-truncated system (\ref{E-tauBnm}) reads as $\dt\overrightarrow{\!\tau} =B_{n,1}\,\ddx\overrightarrow{\!\tau}$.

For small values of $m$, $m=1,2$, the terms $\ddx\tau_{n+m}$ are given by
\begin{eqnarray}\label{E-symm-k1}
 \frac{1}{n+1}\,\ddx\tau_{n+1}\eq\sum\nolimits_{i=1}^n\frac{(-1)^{n-i}}{i}\,\sigma_{n-i+1}\ddx\tau_i
 +\tilde a_1(\overrightarrow{\!\tau},\dx\overrightarrow{\!\tau}),\\
 \label{E-symm-k2}
 \frac{1}{n+2}\,\ddx\tau_{n+2}\eq \sum\nolimits_{i=1}^n\frac{(-1)^{n-i}}{i}
 (\sigma_{1}\/\sigma_{n-i+1}-\sigma_{n-i+2})\ddx\tau_i+\tilde a_2(\overrightarrow{\!\tau},\dx\overrightarrow{\!\tau}).
\end{eqnarray}
By Proposition~\ref{C-BBnm}, the last row of the matrix $B_{n,1}$ or $\frac32(B_{n,1})^2$ consists of the coefficients at $\ddx\tau_i$'s on the right-hand side of (\ref{E-symm-k1}), respectively, of  (\ref{E-symm-k2}), and so on.
\end{remark}

\begin{example}\label{Ex-ex-ric2}\rm
For $f_j=-2\,\delta_{j1}$, (\ref{E-tauBnm}) reduces to the (system of) heat equations
\[
 \dt\tau_i =\ddx\tau_{i},\qquad i\in\NN,
\]
whose solution is known.
Consider more complicated cases.

1. For $f_j=-\delta_{j2}$, (\ref{E-tauBnm}) reduces to the system
\begin{equation}\label{E-tauk-1}
 \dt\tau_i =\frac{i}{i+1}\,\ddx\tau_{i+1},\quad i\in\NN,
\end{equation}
whose $n$-truncated version reads as:
 $\dt\overrightarrow{\!\!\tau}=B_{n,1}\,\ddx\overrightarrow{\!\!\tau}
 +a_2(\overrightarrow{\!\tau},\dx\overrightarrow{\!\tau})$.
 For $n=2$, we have two~PDEs
\[
 \dt\tau_1 = \frac{1}{2}\,\ddx\tau_{2},\quad
 \dt\tau_2 = \frac{2}{3}\,\ddx\tau_{3}
 =(\tau_{2}-\tau_{1}^2)\ddx\tau_{1} +\tau_{1}\ddx\tau_{2},
\]
and the matrix
 $B_{2,1} =\Big(\begin{array}{cc}
           0 & \frac12 \\
          -2\,\sigma_{2} & \sigma_{1} \\
           \end{array}\Big)$.
If the roots $k_1\ne k_2$,
the eigenvectors of $B_{2,1}$ are $v_j=(1,\, 2\,k_j),\ j=1,2$.
For $n=3$, (\ref{E-tauk-1}) reduces to the quasilinear system of three PDEs with the matrix
\[
B_{3,1}=\left(\begin{array}{ccc}
                     0 & \frac12 & 0 \\
                     0 & 0 & \frac23 \\
                     3\,\sigma_{3} & \ \ -\frac{3}{2}\,\sigma_{2} & \ \ \sigma_{1} \\
     \end{array}\right),
\]
whose eigenvalues are $k_{j}$, and the eigenvectors are $v_j=(1, 2\,k_j, 3\,k_j^2)$.

2. For $f_j=-\delta_{j3}$, (\ref{E-tauBnm}) reduces to the system
\begin{equation*}
 \dt\tau_i = \frac{3\,i}{2(i+2)}\,\ddx\tau_{i+2}\quad (i\in\NN)
 \quad\Longleftrightarrow\quad
 \dt\overrightarrow{\!\tau} =-
 \frac{3}2\,(B_{n,1})^{2}\ddx(\overrightarrow{\!\tau})+a_3(\overrightarrow{\!\tau},\dx\overrightarrow{\!\tau}).
\end{equation*}
For $n=3$, the system has the following matrix:
\begin{equation*}
 \frac32 (B_{3,1})^2 =\left(\begin{array}{ccc}
                    0 & 0 & \frac{1}{2} \\
                    3\,\sigma_3 &
                    -\frac{3}{2}\sigma_2 & \sigma_1 \\
                    \ \ \frac{9}{2}\sigma_1\sigma_3 &
                    \ \ \frac{9}{4}(\sigma_3-\sigma_1\sigma_2) &
                    \ \ \frac{3}{2}(\sigma_1^2-\sigma_2) \\
           \end{array}\right),
\end{equation*}
whose eigenvalues are $\frac{3}2\,k_j^{2}$.
This series of examples can be continued as long as one desires.
\end{example}

\end{document}